\documentclass[a4paper]{amsart}

\usepackage{amssymb,amsmath,graphicx,verbatim}
\usepackage{latexsym}
\usepackage{eucal}
\makeatletter
\@addtoreset{equation}{section}\makeatother

\newcommand{\cg}[1]{{}}
\newcommand{\ah}[1]{{}}
\newcommand{\nb}[1]{{}}

\newtheorem{theo}{Theorem}[section]
\newtheorem{lem}[theo]{Lemma}
\newtheorem{prop}[theo]{Proposition}

\newtheorem{rem}[theo]{Remark}

\newcommand{\mc}{\mathcal}
\newcommand{\rr}{\mathbb{R}}
\newcommand{\TT}{\mathbb{T}}
\newcommand{\RR}{\mathbb{R}}
\newcommand{\nn}{\mathbb{N}}

\newcommand{\hh}{\mathbb{H}}
\newcommand{\zz}{\mathbb{Z}}
\newcommand\NN{\mathbb{N}}

\newcommand\Ubar{\overline{U}}

\newcommand\hyp{\operatorname{hyp}}
\def\dbyd#1{\frac{\partial }{\partial #1}}
\newcommand\sech{\operatorname{sech}}

\newcommand\xb{\mathbf{x}}
\newcommand\xib{\mathbf{\xi}}
\newcommand\ad{\operatorname{ad}}

\newcommand{\wf}{\textrm{WF}_h}

\newcommand{\la}{\lambda}
\newcommand{\eps}{\epsilon}

\newcommand{\pl}{\partial}
\newcommand{\x}{\times}

\newcommand{\til}{\widetilde}
\newcommand{\bbar}{\overline}

\newcommand{\supp}{\operatorname{supp}}
\newcommand{\cjd}{\rangle}
\newcommand{\cjg}{\langle}

\newcommand{\demi}{\frac{1}{2}}
\newcommand{\ndemi}{\frac{n}{2}}

\newcommand{\indic}{\operatorname{1\negthinspace l}}

\def\qed{\hfill$\square$}
\newcommand\zbar{\overline{z}}

\newcommand\ang[1]{\langle #1 \rangle}

\begin{document}
\title[Strichartz estimates without loss]{Strichartz estimates without loss on manifolds with hyperbolic trapped geodesics.}
\author{Nicolas Burq}
\author{Colin Guillarmou}
\author{Andrew Hassell}
\begin{abstract}
In \cite{Doi}, Doi proved that the $L^2_tH^{\demi}_x$ local smoothing effect for Schr\"odinger equation  
on a Riemannian manifold does not hold if the geodesic flow has one trapped trajectory. We show in contrast that
Strichartz estimates and $L^1\to L^\infty$ dispersive estimates still hold  without loss for $e^{it\Delta}$ in various situations where the trapped set is hyperbolic and of sufficiently small fractal dimension.
\end{abstract}

\maketitle
The influence of the geometry on the behaviour of solutions of linear or non linear partial differential equations has been widely studied recently, and especially in the context of wave or Schr\"odinger equations. In particular, the understanding of the smoothing effect for the Schr\"odinger flow and Strichartz type estimates has been related to the global behaviour of the geodesic flow on the manifold (see for example the works by Doi~\cite{Doi} and Burq~\cite{Burq1}). Let us recall that for the Laplacian $\Delta$ on a $d$-dimensional non-compact Riemannian manifold $(M,g)$, the \emph{local smoothing effect}
for bounded time $t\in[0,T]$ and Schr\"odinger waves  $u = e^{it\Delta} u_0: M \times \RR \to \mathbb{C}$
 is the estimate 
\[||\chi e^{it\Delta}u_0||_{L^2((0,T);H^{1/2}(M))}\leq C_T||u_0||_{L^2(M)}, \quad \forall u_0\in L^2(M)\]
where $C_T>0$ is a constant depending a priori on $T$ and $\chi$ is a compactly supported smooth function (the assumption on $\chi$ can of course be weakened in many cases, e.g for $M=\rr^d$) \cite{CS}.  
In other words, although the solution is only $L^2$ in space uniformly in time, it is actually half a derivative better (locally) in an $L^2$-in-time sense. 
For its description in geometric settings, the picture now is fairly complete:  the so called ``nontrapping condition'' stating roughly that every geodesic maximally extended goes to infinity, is known  to be essentially necessary and sufficient (modulo reasonable conditions near infinity) \cite{Burq1}. 

Another tool for analyzing non-linear Schr\"odinger equations is the family of so-called \emph{Strichartz estimates} introduced by \cite{St}: for Schr\"odinger waves on Euclidean space $\RR^d$ with initial data $u_0$,
\begin{equation}\label{firstST}
||e^{it\Delta}u_0||_{L^p((0,T);L^q(\RR^d))}\leq C_T||u_0||_{L^2(\RR^d)}\quad {\rm if}\quad p, q \geq 2, \quad \frac{2}{p}+\frac{d}{q}=\frac{d}{2}, \quad (p,q)\not=(2,\infty).
\end{equation}  
If $\sup_{T\in(0,\infty)}C_T<\infty$, we will say that a global-in-time Strichartz estimate holds. Such a global-in-time estimate has been proved by Strichartz for the flat Laplacian on $\rr^d$ while the local-in-time estimate is known in several geometric situations where the manifold
is non-trapping (asymptotically Euclidean, conic or hyperbolic); see \cite{BT,B,HTW,ST}.
On the other hand it is clear that such a global-in-time estimate cannot hold on compact manifolds, for
it suffices to consider the function $u_0=1$.
The situation is similar for the non-compact case in the presence of  elliptic (stable) non-degenerate periodic orbits of the geodesic flow: as remarked by M. Zworski, the quasi-modes constructed by  Babi\v c~\cite{Bab} and Py\v skina~\cite{Py}  (see also the work by Ralston~\cite{Ra}) show that for Schr\"odinger solutions, some loss must occur as far as Strichartz (or smoothing) estimates are concerned; and moreover, that no Strichartz estimates can be true globally in time in the presence of such orbits.
On the other hand Burq-G\'erard-Tzvetkov \cite{BGT} proved that \eqref{firstST} holds on compact manifolds for finite time
if one replaces $||u_0||_{L^2(M)}$ by $||u_0||_{H^{1/p}(M)}$, meaning that a Strichartz estimate is satisfied if one accepts 
some loss of derivatives. It is however certainly not optimal in general since Bourgain \cite{Bou} proved that for the flat torus $(\rr/2\pi\zz)^2$, 
the Strichartz estimate for $p=q=4$ holds with $\eps$ loss of derivatives for any $\eps>0$.  
Another striking example has been given by Takaoka and Tzvetkov \cite{TT} by adapting the ideas of Bourgain, 
namely the case of the two dimensional infinite flat cylinder 
$S^1\x \rr$ where \eqref{firstST} holds (with no loss of derivatives) if $p=q=4$; note that this manifold is trapping.  
An example with a repulsive potential $V(x_1,x_2)=x_1^2-x_2^2$ has also been studied by Carles \cite{Car},
who proved that global-in-time Strichartz estimates with no loss hold in this case.  To summarize, it is not really understood when (i.e. under what geometric conditions) a loss in Strichartz estimates must occur, and if it does, how large that loss must be.

The purpose of this article is precisely to give some examples of Riemannian manifolds where trapping does occur (and consequently loss is unavoidable for the smoothing effect), but nevertheless, since the dynamics are \emph{hyperbolic} near the trapped set, 
we are able to prove (local-in-time) Strichartz estimates without loss for Schr\"odinger solutions.\\ 

The first example, which we treat in Section 1,  
is the case of a convex co-compact hyperbolic manifold of dimension $d=n+1$, with a limit set of Hausdorff dimension $\delta<n/2$. The simplest example of such a manifold is the two dimensional infinite hyperbolic cylinder with one single
trapped geodesic. In this case, the calculations are quite explicit, representing the Schr\"odinger kernel as an average over the group of the Schr\"odinger kernel on the hyperbolic space $\hh^{n+1}$,  and we are able to prove that not only Strichartz estimates, but the stronger $L^1 \to L^\infty$ dispersive estimates  hold for the Schr\"odinger group.
\begin{theo}
Let $X$ be an $(n+1)$-dimensional convex co-compact hyperbolic manifold such that its limit set has Hausdorff dimension $\delta<n/2$.
Then the following dispersive and Strichartz estimates without loss hold:
\[\begin{gathered}
||e^{it\Delta_X}||_{L^1(X)\to L^\infty(X)}\leq 
\begin{cases} C |t|^{-(n+1)/2}, \textrm{ for }|t|\leq 1 \\
C |t|^{-3/2}, \textrm{ for }|t|>1,\end{cases}\\
||e^{it\Delta_X}u_0||_{L^p(\rr;L^q(X))}\leq C||u_0||_{L^2(X)}
\end{gathered}\]
for all $(p,q)$ such that $(1/p,1/q)\in T_n$ where  
\begin{equation}\label{admissible}
T_n:=\left\{\Big(\frac{1}{p},\frac{1}{q}\Big)\in \Big(0,\demi\Big]\x\Big(0,\demi\Big); \frac{2}{p}\geq \frac{n+1}{2}-\frac{n+1}{q} \right\}\cup \left\{\Big(0,\demi\Big)\right\}.
\end{equation}
\end{theo}
These manifolds are non-compact, infinite volume, with finitely many ends of funnel type, they have constant curvature $-1$  and possess \emph{infinitely} many closed geodesics; it is remarkable that despite this last fact, a sharp dispersive estimate holds for all time. 
We also remark that Strichartz estimates for the same range of $(p,q)$ have been recently 
shown by Anker-Pierfelice \cite{AP} for the model non-trapping case $\hh^{n+1}$ (see also \cite{Ba,IoS} for the estimate
\eqref{firstST} in that setting). The
triangle of admissibility for the Strichartz exponent $(p,q)$ is a consequence of the exponential decay of the 
integral kernel of the Schr\"odinger operator at infinity. Notice that in the asymptotically hyperbolic setting, J-M. Bouclet \cite{B} proved Strichartz estimates without loss of derivatives for bounded times and with admissibility exponents satisfying \eqref{firstST}  for \emph{non-trapping} such manifolds (in this case the sectional curvature is not assumed constant, but rather tending to $-1$ at infinity).\\  

The second example is on the manifold $Z$ given by the connected sum of two copies of Euclidean $\RR^2$. This we provide with a Riemannian metric $g$ by gluing two copies of the Euclidean metric on $\RR^2$ with the metric on the 2-dimensional hyperbolic cylinder. Essentially because $Z$ is formed from pieces all of which satisfy Strichartz estimates without loss, the same is true for $(Z,g)$. Actually we need to use local smoothing estimates to control error terms in the transition region, but since this region is disjoint from the single trapped orbit, there are no losses in such local smoothing estimates. This example is given in Section 2; the main result is Theorem~\ref{example}. \\

Our last family of examples, in Section 3, is a generalization of that in Section 2 to higher dimensions and more complicated trapped sets. It is similar to the class of manifolds studied recently by Nonnenmacher and Zworski~\cite{NZ}: we  
consider asymptotically Euclidean (or more generally asymptotically conic) manifolds, 
the curvature of which is assumed to be negative in a geodesically convex compact part that includes
the (projection of) the trapped set, and such that the trapped set is small enough in the sense that 
the topological pressure $P(1/2)$ of the trapped set evaluated at $1/2$ is negative (see the definition in Section~\ref{tp}). 
This  is a dynamical condition which generalizes the condition $\delta<n/2$ above to this more general setting and
roughly speaking means that the trapped set is filamentary with small fractal dimension. For instance, for surfaces (dimension $d=2$) this means
that the trapped set (as a subset of the cosphere bundle $S^*M$) has Hausdorff dimension less than $2$.
More precisely our result (which includes the example in Section 2 as a special case) is 
\begin{theo}\label{zzzz} Let $(M, g)$ satisfy assumptions (A1) --- (A4) defined in \eqref{negpress}. Then Strichartz estimates without loss hold for $M$:  there exists $C>0$ such that
\begin{equation}
||e^{it\Delta}u_0||_{L^p((0,1),L^q(M))}\leq C||u_0||_{L^2(M)}
\end{equation}
for all $u_0\in L^2(M)$ and $(p,q)$ satisfying \eqref{firstST} and $p > 2$.
\end{theo}
Note that Christianson \cite{Chr} and Datchev \cite{Da} showed that Strichartz estimates hold with an $\eps$ loss of derivatives for all $\eps>0$ in that setting.  Our method is based on the use of the local smoothing effect with log loss, which follows from the
resolvent estimate of \cite{NZ} (see also \cite{Da})
\[||\chi(\Delta-\la+i0)^{-1}\chi||_{L^2\to L^2}\leq C\frac{\log(\la)}{\la^\demi}, \quad \textrm{ for }\chi\in C_0^\infty(M)\]
together with a sharp dispersive estimate on the logarithmically \emph{extended} time interval $t\in (0,h\log(h))$ for the frequency localized operator $e^{it\Delta}\varphi(h^2\Delta)$ 
where $\varphi\in C_0^\infty((0,\infty))$ and $h\in(0,h_0)$ is small. Roughly speaking this logarithmic extension of the time interval of validity of the dispersive estimate allows one to recover
the $\log$ loss in the local smoothing estimate. This dispersive estimate is inspired by the works of Anantharaman \cite{An}, Anantharaman-Nonnenmacher~\cite{AnNo} and
Nonnenmacher-Zworski \cite{NZ}. In particular, the technique for proving the dispersive estimate for logarithmically extended time originates in \cite[Theorem 1.3.3]{An}, while the idea of combining the exponential decay provided by this theorem with the topological pressure assumption (see (A4) in Section 3) is due to \cite{NZ}.\\ 

\textbf{Acknowledgement}. We thank S. Nonnenmacher, N. Anantharaman and F. Planchon for helpful 
discussions and references. N.B. is supported by ANR grant ANR-07-BLAN-0250. C.G. is supported by ANR grant ANR-09-JCJC-0099-01 and thanks the Mathematical Sciences Institute of 
ANU Canberra where part of this work was done. A.H. is supported by Australian Research Council Discovery Grant DP0771826 and thanks the mathematics department at Universit\'e Paris 11 for its hospitality. We are finally grateful to the referee for his 
careful reading.


\section{Hyperbolic manifolds}\label{hypman}
A \emph{convex co-compact subgroup} $\Gamma\subset {\rm SO}(n+1,1)$ is a discrete group of orientation preserving isometries of hyperbolic space $\hh^{n+1}$, consisting of hyperbolic isometries and such that the quotient $X:=\Gamma\backslash\hh^{n+1}$ 
has finite geometry and infinite volume. If one considers the ball model $B^{n+1}$ of $\hh^{n+1}$, a 
hyperbolic isometry is an isometry of $\hh^{n+1}$ which fixes exactly two points on $\bbar{B^{n+1}}$, and these points are
on the boundary $S^n=\pl\bbar{\hh^{n+1}}$.
The manifold $X:=\Gamma\backslash\hh^{n+1}$ is said to be \emph{convex co-compact hyperbolic}; it is a smooth complete hyperbolic manifold which admits a natural conformal compactification $\bar{X}$ and the hyperbolic metric $g$ on $X$ is of the form
$g=\bar{g}/x^{2}$ where $x$ is a smooth boundary defining function of $\bar{X}$ and $\bar{g}$ a smooth metric on $\bar{X}$. 
The set of closed geodesics is in correspondence with the classes of conjugacy of the group $\Gamma$. The limit set of $\Gamma$ is the set of accumulation points on the sphere $S^{n}=\pl\bbar{\hh^{n+1}}$ of the orbit $\Gamma.m$ where $m\in\hh^{n+1}$ is any point. It has a Hausdorff dimension given by $\delta\in [0,n)$,  and the trapped set of the geodesic flow on the unit tangent bundle $SX$ has Hausdorff dimension $2\delta+1$; see \cite{Su,Zw}.

The simplest example is $\Gamma = \mathbb{Z}$ acting by powers of a fixed dilation $D$ on the upper half space model of $\hh^{n+1}$. Then the limit set consists of two points $\{ 0, \infty \}$, $\delta = 0$, and $\hh^{n+1}/\Gamma$ is the $(n+1)$-dimensional hyperbolic cylinder. 

We prove the following.

\begin{theo}\label{stric}
Let $X$ be an $(n+1)$-dimensional convex co-compact hyperbolic manifold such that its limit set has Hausdorff dimension $\delta<n/2$.
Then $e^{it\Delta_X}$ has a smooth Schwartz kernel for $t \neq 0$, and there is a constant $C$ such that the following dispersive estimate holds for all $t \neq 0$:
\begin{equation}\label{disp}
|| e^{it\Delta_X}||_{L^1\to L^\infty}\leq 
\begin{cases} C |t|^{-(n+1)/2}, \textrm{ for }|t|\leq 1 \\
C |t|^{-3/2}, \textrm{ for }|t|>1.\end{cases}
\end{equation}
Moreover the following global-in-time Strichartz estimates hold:
\begin{equation}
||e^{it\Delta_X}u_0||_{L^p(\rr;L^q(X))}\leq C||u_0||_{L^2(X)},
\label{gse}\end{equation}
for all $(p,q)$ such that $(1/p,1/q)\in T_n$ where  $T_n$ is given by \eqref{admissible}.
\end{theo}

\begin{proof}
The integral kernel of the Schr\"odinger operator $e^{it\Delta_{\hh^{n+1}}}$ on  
hyperbolic space has been computed by V. Banica \cite{Ba}. It is a function of the hyperbolic distance 
\begin{equation}\label{kernel}
\begin{gathered}
K(t;\rho(z,z'))=c|t|^{-\demi}e^{-itn^2/4}(\sinh(\rho)^{-1}\pl_\rho)^{\ndemi}e^{i\rho^2/4t}, \quad n \textrm{ even}\\
K(t;\rho(z,z'))=c|t|^{-\frac{3}{2}}e^{-itn^2/4}(\sinh(\rho)^{-1}\pl_\rho)^{\frac{n-1}{2}}\int_{\rho}^\infty \frac{e^{is^2/4t}s}{\sqrt{\cosh s-\cosh \rho}}ds, \quad n\textrm{ odd}
\end{gathered}
\end{equation}
where $\rho=\rho(z,z'):=d_{\hh^{n+1}}(z,z')$.
In both cases we remark, like for the heat kernel, that the kernel $K(t;.,.)$ is smooth 
on $\hh^{n+1}\x\hh^{n+1}$ for $t\not=0$; this is clear when $n+1$ is odd, and needs a bit more 
analysis when $n+1$ is even.
From this expression we obtain an upper bound for $|K(t;z,z')|$  (see \cite[Prop 4.1 and Sec. 4.2]{Ba}) for $t\not=0$  of the form 
\begin{equation}\label{boundk} 
\begin{cases} 
C |t|^{-(n+1)/2}\Big(\frac{\rho}{\sinh\rho}\Big)^{\ndemi} , \textrm{ for }|t|\leq 1\\
C|t|^{-3/2}
\Big(\frac{\rho}{\sinh\rho}\Big)^{\ndemi}, \ \textrm{ for }|t|>1
\end{cases}\end{equation}
for some constant $C>0$. 
Using the inequality $(\rho/\sinh \rho) \leq (1 + \rho) e^{-\rho}$, and since 
$$i\pl_tK(t;z,z')=-\Delta_zK(t;z,z')=-\Delta_{z'}K(t;z,z'),$$
 one can deduce that for $t\not= 0$ bounded 
\begin{multline}\label{delta}
|\Delta^j_zK(t;z,z')| + |\Delta^j_{z'}K(t;z,z')|\leq 
\begin{cases} 
C'|t|^{-(n+1)/2-2j}(1+\rho)^{\ndemi+2j} e^{-\ndemi \rho} ,  \textrm{ for }|t|\leq 1\\
C'|t|^{-3/2}(1 + \rho)^{\ndemi+2j} e^{-\ndemi \rho}, \textrm{ for }|t|>1
\end{cases}
\end{multline}
and in particular $K(t, z, z')$ is smooth in $z,z'$.

To proceed we use the celebrated result of 
Patterson and Sullivan \cite{Pa,Su} that the dimension of the limit set $\delta$ is the exponent of convergence of the Poincar\'e series
\[P_s(z,z'):=\sum_{\gamma\in \Gamma}e^{-s\rho(z,\gamma.z')}, \quad z,z'\in\hh^{n+1}.\] 
\begin{lem}
Let $\mc{F}\in\hh^{n+1}$ be a fundamental domain of the convex co-compact group $\Gamma$ 
and let $x$ be a boundary defining function of the compactification $\bar{X}$ of $X=\Gamma\backslash\hh^{d+1}$, which we also
view as a function on $\mc{F}$. For each $\gamma\in\Gamma$, define by $\ell_\gamma$ the translation length of $\gamma$. 
Then there exists $R>0$ such that for all $\eps>0$, there is $C_\eps>0$ such that 
for $s>\delta+\eps$ and all $z,z'\in \mc{F}$,
\begin{equation}\label{boundp}
\sum_{\gamma\in\Gamma, \ell_\gamma>R}e^{-s\rho(z,\gamma z')}\leq C_\eps (x(z) x(z'))^s.
\end{equation} 
\end{lem}
\begin{proof} In \cite[Lemma 5.2]{GMP}, it is shown that there exist constants $C_1,C_2>0$ such that 
for all $\gamma\in\Gamma$ such that $\ell_\gamma>C_1$ 
\[    e^{-\rho(z,\gamma z')}\leq C_2e^{-\ell_\gamma}x(z)x(z').\] 
Now it suffices to sum after raising to the power $s$ and to use the fact 
that $\sum_{\gamma\in\Gamma}e^{-s\ell_\gamma}<C_\eps$ for some $C_\eps$ if $s>\delta+\eps$.
\end{proof}

Combining (\ref{boundp}) and (\ref{boundk}), we deduce that for $z,z'\in\mc{F}$, the series 
\[K_X(t;z,z'):=\sum_{\gamma\in\Gamma}K(t;z,\gamma.z')\] 
converges uniformly and for all $s<n/2$ there exists $C_s>0$ such that 
for all $z,z'\in\mc{F}$
\begin{equation}\label{pointwise}
|K_X(t;z,z')|\leq 
\begin{cases} C_s|t|^{-(n+1)/2}\Big(\sum_{\gamma\in\Gamma,\ell_\gamma\leq R}
e^{-s\rho(z,\gamma z')}+x(z)^{s}x(z')^s\Big), \textrm{ for }|t|\leq 1 \\
C_s|t|^{-3/2}\Big(\sum_{\gamma\in\Gamma,\ell_\gamma\leq R}e^{-s\rho(z,\gamma z')}+x(z)^{s}x(z')^s\Big), 
\textrm{ for }|t|> 1 \end{cases}
\end{equation}
where $R$ is the constant in \eqref{boundp}. This leads directly to the dispersive estimate  
\begin{equation}\begin{gathered}
\sup_{z,z'\in \mc{F}}|K_X(t;z,z')|\leq 
\begin{cases} C'|t|^{-(n+1)/2}, \textrm{ for } |t|\leq 1 \\
 C'|t|^{-3/2} , \textrm{ for } |t|>1 \end{cases}.
\end{gathered}\end{equation}
for some constants $C'$.
Moreover, using (\ref{delta}), the same argument shows that 
the series $K_X(t;z,z')$ is smooth in $z,z'$ for $t\not=0$.
Let $\mc{F}$ be a fundamental domain of $\Gamma$. 
For any $u_0\in C_0^{\infty}(X)$ the function $u(t):=\int_{\mc{F}}K_X(t;z,z')u_0(z')dz'$ is smooth on $\hh^{n+1}$ and
satisfies $u(t, \gamma z) = u(t, z)$ for any $\gamma\in\Gamma$, 
thus $u(t)$ is smooth on $X$. Moreover it solves the Schr\"odinger equation on $X$ with initial data $u(0)=u_0$, so $u(t)=e^{it\Delta_X}u_0$.
This implies that $K_X$ is the Schwartz kernel of $e^{it\Delta_X}$ on $X$.\\
 
 We next prove the global-in-time Strichartz estimates \eqref{gse} (notice that these estimates for a \emph{finite} time interval follow immediately from the small-time dispersive estimate \eqref{disp} and from Keel-Tao \cite{KT}, following the method of Anker-Pierfelice \cite{AP}). 
Let $\gamma\in\Gamma$ be such that $\ell_\gamma\leq R$ where $R$ is the constant in \eqref{boundp}, then  define
$K_\gamma^{1}(t)$ to be the operator acting on $\mc{F}$ with $L^\infty$ kernel $\indic_{\mc{F}}(z)K(t;z,\gamma z')\indic_{\mc{F}}(z')$. Since $\gamma$ is an isometry of $\hh^{n+1}$, 
this operator can also be written as 
$f\to \indic_{\mc{F}}e^{it\Delta_{\hh^{n+1}}}\gamma_*(\indic_{\mc{F}}f)$.  
Then from Theorem 3.4 of Anker-Pierfelice \cite{AP} and the fact that push-forward $\gamma_*$ is an isometry on any
$L^{r'}(\hh^{n+1})$, we get the estimate 
\begin{equation}\label{lqlq'}
||K_\gamma^{1}(t)u_0||_{L^q(\hh^{n+1})}\leq C||\indic_{\mc{F}} u_0||_{L^{r'}(\hh^{n+1})}
\x \left\{\begin{array}{ll}
|t|^{-(n+1)\max(\demi-\frac{1}{q},\demi-\frac{1}{r})} & \textrm{ if }|t|\leq 1\\
|t|^{-\frac{3}{2}} & \textrm{ if }|t|>1
\end{array}\right.
\end{equation}
for all $2<q,r\leq \infty$ and $1/r'+1/r=1$, so the same estimate holds for 
$K^1(t):=\sum_{\gamma\in\Gamma, \ell_\gamma\leq R}K^1_\gamma(t)$.
Now consider the operator $K^2(t)$ acting on $\mc{F}$ whose $L^\infty$ 
kernel is $K^2(t;z,z'):=K_X(t;z,z')-K^1(t;z,z')$. From (\ref{boundp}) and (\ref{boundk}), this kernel is bounded (for all $s<n/2$) by 
\[|K^2(t;z,z')|\leq 
\begin{cases} C_s|t|^{-(n+1)/2}x(z)^{s}x(z')^s, \textrm{ for } |t|\leq 1 \\
C_s|t|^{-3/2}x(z)^{s}x(z')^s , \textrm{ for } |t|>1 \end{cases}.\]
Since the hyperbolic metric on $\mc{F}$ induces a measure of the form $x^{-n-1}\mu$ for some
bounded measure $\mu$ on $\mc{F}$ we see that the function $x^{s}$ is in $L^{\alpha}(\mc{F},{\rm dv}_{\hh^{n+1}})$ for 
all $\alpha>n/s$, and hence deduce directly that $K^2(t)$ satisfies
\[||K^2(t)u_0||_{L^q(\mc{F},{\rm dv}_{\hh^{n+1}})}\leq C||u_0||_{L^{r'}(\mc{F},{\rm dv}_{\hh^{n+1}})}
\x \left\{\begin{array}{ll}
|t|^{-\frac{n+1}{2}} & \textrm{ if }|t|\leq 1\\
|t|^{-\frac{3}{2}} & \textrm{ if }|t|>1
\end{array}\right.\] 
for all $q,r\in(2,\infty]$ and $r'$ the conjugate exponent of $r$,
so the same estimate holds for $K_X(t)$ when combining with \eqref{lqlq'}.
Then it suffices to conclude using the standard $TT^*$ argument exactly as in
the proof of Theorem 3.6 of Anker-Pierfelice \cite{AP} 
and we obtained the claimed Strichartz estimate.
\end{proof}

\begin{rem} 
For hyperbolic quotients with dimension of limit set $\delta>n/2$, the positive number $\delta(n-\delta)$
is an $L^2$ eigenvalue with multiplicity one and smooth eigenvector $\psi_\delta$. It follows easily from Section 2 of \cite{Pe} (or the general result of Mazzeo-Melrose \cite{MM} about the structure of the resolvent)
that $\psi_\delta\in L^p(X)$ for all $p>n/\delta$, thus in particular for all $p\geq 2$. This implies 
that for $q\geq 2$, $2<p<\infty$ and all $\chi\in L^\infty(X)$
\[\|\chi e^{it\Delta_X}\psi_\delta\|_{L^q(X)}=\|\chi \psi_\delta\|_{L^q(X)} \notin L_t^p((0,\infty)),\]
so global-in-time Strichartz estimates cannot hold when $\delta>n/2$, even with a space cut-off.
\end{rem}

\section{Connected sum of two copies of $\rr^2$}

In this section, we give an example of a Riemannian manifold $(Z,g)$ which is topologically the connected sum of two copies of $\RR^2$, and is geometrically Euclidean near infinity, and hyperbolic near the `waist' (and hence with a single trapped ray),  for which Strichartz estimates without loss are valid. The idea is simple; since Strichartz estimates without loss are valid on flat $\RR^2$, and on the hyperbolic cylinder (thanks to Theorem~\ref{stric}), then they should also be valid on a space obtained by gluing pieces of these manifolds together, provided that no additional trapping is created by the gluing procedure. 

Let us consider an asymptotically Euclidean manifold $(Z, g)$ which is the connected sum of two copies of $\RR^2$, joined by a neck which has a neighbourhood $U$  isometric to a neighbourhood $U'$ of the short closed geodesic, or `waist', on the hyperbolic two-cylinder $C^2$. We denote this short closed geodesic by $\gamma$, whether on $Z$ or on $C^2$. We can write down an explicit metric $g$ for such a manifold, on $\RR \times S^1$, in the form $dr^2 + f(r)^2 d\theta^2$, where $d\theta^2$ is the metric on $S^1$ of length $2\pi$, and where $f(r) = \cosh r$ for small $r$, say $r \leq 3\eta$ for some small $\eta > 0$, and is equal to $|r| + a$ for large $|r|$, say $|r| \geq R$ (where $a$ is a constant). We also choose $f$ so that $f'(r)$ has the same sign as $r$; it is easy to see that this is compatible with the condition that $f(r) = \cosh r$ for small $|r|$ and $|r| + a$ for large $|r|$. 
The equations of motion for geodesic flow then give 
$\ddot r = 2 f'(r) f(r)  \dot \theta^2$, which has the same sign as $r$, and it is straightforward to deduce from this that there can be no trapped geodesic other than the waist $\gamma$ at $r=0$. 
For any such manifold $(Z,g)$ we have

\begin{theo}\label{example} For any finite $T$ there is a constant $C_T$ such that 
\begin{equation}\label{Strichartz-example}
\| e^{it\Delta_Z} u_0 \|_{L^p([0,T]; L^q(Z))} \leq C_T \| u_0 \|_{L^2(Z)}
\end{equation}
for all $(p,q)$ satisfying \eqref{firstST} with $d=2$ and all $u_0 \in L^2(Z)$. 
\end{theo}

Before giving the proof we introduce some further notation and definitions. We will compare $Z$ to the hyperbolic cylinder $C^2$ and to the auxiliary Riemannian manifold $(\tilde Z = \RR^2, \tilde g)$, given in standard polar coordinates $(r, \theta)$ on $\RR^2$ by $\tilde g = dr^2 + \tilde f(r)^2 d\theta$, where $\tilde f(r) = f(r)$ for $r \geq \eta$, $\tilde f'(r)>0$, and is equal to $r$ for small $r$. Reasoning as above, we see that the metric $\tilde g$ is nontrapping. 
We will take a Schr\"odinger wave $u$ on $Z$ and decompose it so that one piece lives on $C^2$ and the other lives on $\tilde Z$, and we will deduce Strichartz without loss on $Z$ from the fact that Strichartz without loss holds for both $C^2$ and $\tilde Z$.

\begin{proof} Let $U \subset Z$ be a neighbourhood of $\{ r=0 \}$, say $\{ |r| < 2\eta \}$, thus containing the projection of the trapped set of $Z$. By construction, the metric is exactly hyperbolic in a neighbourhood of  $\Ubar$. 
We decompose $u = u_i + u_e$, where $u_i = \chi u$ is supported in $U$ and $u_e = (1 - \chi) u$ is supported where the metric $g$ is identical to $\tilde g$. (Thus $\nabla \chi$ is supported where $\eta \leq |r| \leq 2\eta$.) We prove the estimate \eqref{Strichartz-example} separately for $u_i$ and $u_e$. As stated above, the idea is to regard $u_e$ as solving a PDE on $\tilde Z$ and to regard $u_i$ as solving a PDE on $C^2$. 

We first prove a local smoothing result for $Z$, $\tilde Z$ and $C^2$. This is essentially standard, but we give the details for the reader's convenience (and in keeping with the expository character of this section). For applications in the following section, we 
give a result in any dimension.

\begin{lem}\label{T} (i) Suppose that $X$ is a $d$-dimensional manifold with Euclidean ends and with trapped set $K \subset T^*X$. Suppose that $u$ solves the Schr\"odinger equation on $(X,g)$ with initial condition $u_0$, and suppose that $\phi \in C_c^\infty(M)$ is supported away from the projection of the trapped set $\pi(K)$. Then 
$$
\| \phi u \|_{L^2_t([0,T]; H^{\demi}(X))} \leq C \| u_0 \|_{L^2(X)}.
$$

(ii) Suppose that $v$ solves the Schr\"odinger equation on $(C^2, g_{\hyp})$ with initial condition $v_0$, and suppose that $\phi \in C_c^\infty(C^2)$ is supported away from the closed geodesic $\gamma$. Then 
$$
\|\phi v \|_{L^2_t([0,T]; H^{\demi}(C^2))} \leq C \| \tilde v_0 \|_{L^2(C^2)}.
$$

\end{lem}

\begin{proof} This result can be deduced from the resolvent estimate of Cardoso-Vodev \cite{CV}, but for completeness we give a proof via a positive commutator argument. We construct a zeroth order pseudodifferential operator $A$ on $X$ such that $i[\Delta, A]$ has a nonnegative principal symbol which is elliptic on the support of $\phi$. Then we use the identity
\begin{equation}
\ang{Au(\cdot, T), u(\cdot, T)} - \ang{Au(\cdot, 0), u(\cdot, 0)}
= \int_0^T \ang{i[\Delta, A] u, u} \, dt 
\label{identity}\end{equation}
valid for any Schr\"odinger wave $u$. Since $i[\Delta, A]$ is order one and elliptic on the support of $\phi$, the right hand side is equal to $c \| \phi u \|_{L^2([0,T]; H^{1/2})}^2$ plus terms which are essentially positive, while the left hand side is bounded by $C \| u_0 \|_{L^2}^2$, giving the estimate. 

Let us set $A = A_1 + A_2$, where $A_1$ is supported in the region where $X$ is Euclidean and $A_2$ is properly supported.  We take $R > 0$ large enough so that each end of $X$ has a neighbourhood isometric to $\RR^d \setminus \overline{B(0,R)}$ and use Euclidean coordinates $\xb = (x_1, x_2, \dots x_d)$ on this neighbourhood with dual cotangent coordinates $\xib$. We write $r = |\xb|$ and take $A_1$ to have principal symbol 
\begin{equation}
a_1 = \zeta^2(r) \ang{\xib}^{-1} r^{-1} \xb \cdot \xib \Big( 1 - r^{-\epsilon} \Big).
\label{a1}\end{equation}
Here $\zeta(t)$ is chosen to be $0$ for $t < R$ and $1$ for $t \geq 2R$ and to be nondecreasing, where $R$ is sufficiently large that $1 - R^{-\epsilon} > \demi$, say. 
We understand this to mean that $a_1$ is defined as above on \emph{each} end of $X$. 
Explicitly, we could take $A_1 = \zeta(r) r^{-1} (1 + \Delta)^{-1/2} \demi(\xb \cdot D_\xb + D_\xb \cdot \xb) \zeta(r) (1 - r^{-\epsilon})$, where here $\Delta$ denotes the flat Laplacian on $\RR^d$; notice that $(1 + \Delta)^{-1/2}$ makes sense since it is both pre- and post-multiplied by $\zeta(r)$ which is supported where the metric is Euclidean. 

Then the derivative $a_1$ along the Hamilton vector field of $\sigma(\Delta_X)$, namely the geodesic flow $2\xib \cdot \partial_\xb$, is 
$$\begin{gathered}
2\zeta^2(r) \ang{\xib}^{-1} r^{-3}  \big( 1 - r^{-\epsilon} \big) \Big( r^2 |\xib|^2 - (\xb \cdot \xib)^2 \Big) \\
+4 \zeta(r) \zeta'(r) \ang{\xib}^{-1} r^{-2}  \big( 1 - r^{-\epsilon} \big) (\xb \cdot \xib)^2 \\
+ 2\zeta^2(r) \ang{\xib}^{-1} r^{-3} \epsilon r^{-\epsilon} (\xb \cdot \xib)^2 .
\end{gathered}$$
We see that this is nonnegative everywhere, and bounded below by $C \ang{\xib}r^{-1 - \epsilon}$ for $r \geq 2R$. 

Now we define a symbol $a_2$ which will be supported in the region $r \leq 4R$. First we introduce some notation: for $\tilde R \geq R$, let
$E_{\tilde R}$ denote the union of ends $\RR^d \setminus B(0, \tilde R)$, and let $U_{\tilde R} \subset T^*X$ denote $\pi^{-1}(X \setminus E_{\tilde R})$. We choose conic neighbourhoods $U_<$ and $U$ of $K$ such that $\overline{U_<} \subset U$ and $\pi(U)$ is disjoint from $\supp \phi$. For any $p \in \pi^{-1}(\supp \phi)$ we let $\beta = \beta(p)$ denote the maximally extended geodesic through $p$, and we denote by $\beta_+$, resp. $\beta_-$, the forward, resp. backward geodesic ray starting at $p$. Standard topological arguments show that one can choose $U_<$ so that for every $p \in \pi^{-1}(\supp \phi)$, at least one of $\beta_+$ or $\beta_-$ does not meet $U_<$, which we will now assume. 

Now choose an arbitrary $p \in \pi^{-1}(\supp \phi)$ and consider the geodesic $\beta(p)$.  Because of the way we chose $U_<$, either $\beta_-$ or $\beta_+$ does not intersect $\overline{U_<}$. Suppose, for the sake of definiteness that $\beta_+$ does not meet $\overline{U_<}$ (the argument for $\beta_- \cap \overline{U_<} = \emptyset$ is similar). Let $V$ be a conic neighbourhood of $\beta_+$. We may construct a symbol of order $0$ that is 
\begin{itemize}
\item supported on $V$, 
\item non-decreasing with respect to geodesic flow on $U_{2R}$,
\item strictly increasing with respect to geodesic flow on $V_< \cap (U_{2R} \setminus U) \cap \{ |\xi| \geq 1/2 \}$, where $V_<$  is a conic neighbourhood of $\beta_+$ such that $\overline{V_<} \subset V$, and 
\item vanishing outside $U_{4R}$.
\end{itemize}
To do this, we let $t$ be an arc-length parameter along $\beta$ with $\beta(0) = p$ and $\beta_+ = \{\beta(t) \mid t \geq 0 \}$, and let $t_1 = \sup \{ t \mid \beta(t) \in \overline{U_<} \}$, which is negative by assumption. Also let $t_2 = \sup \{ t \mid \beta(t) \in U_{2R} \}$, and $t_4 = \sup \{ t \mid \beta(t) \in U_{4R} \}$, both of which are positive by assumption. 
We choose a function along $\beta$ that is $0$ for $t \leq t_1$, strictly increasing for $t_1 < t < t_2$ and zero for $t \geq t_4$. This can be extended to a symbol of order 0 supported in $V$.

Using compactness, we can select a finite number of conic neighbourhoods $V_<$ as above, covering $\pi^{-1} \supp \phi \setminus \{ 0 \}$. Summing the corresponding symbols defined above, we obtain a symbol $a_2$ supported in $U_{4R}$ such that the Hamilton vector field of $\Delta_X$ applied to $a_2$ is positive and elliptic on $\pi^{-1} \supp \phi$. Let $A_2$ be a properly supported pseudodifferential operator with symbol $a_2$. Let $A$ be the sum of $A_2$ and a sufficiently large multiple of $A_1$. Then $i[\Delta,A]$ has nonnegative symbol, and (if the symbol $a_2$ is specified appropriately, i.e. so that $\{ \sigma(\Delta), a_2 \}$ is a sum of squares of symbols, which is always possible) may be expressed in the form $\sum_i B_i^* B_i + B_0$, where the $B_i$ are order $1/2$ and $\sum_i B_i$ is elliptic on $\pi^{-1} \supp \phi$ and $B_0$ is order $0$. Then substituting $i[\Delta,A] = \sum_i B_i^* B_i + B_0$, and using the sharp G{\aa}rding inequality in the form $C \sum_i \| B_i u \|_2^2 \geq  \| \phi u \|_{H^{1/2}(X)}^2 - C' \| u \|_{L^2(X)}^2$ which is valid for sufficiently large $C$,   we deduce that
$$
\int_0^T \| \phi u \|_{H^{\demi}(M)}^2 \, dt \leq C \| u_0 \|_{L^2(M)}^2,
$$
proving (i).

The proof of (ii) is very similar in spirit. Again we construct a pseudodifferential operator $A$ with the property that $i[\Delta, A]$ has a nonnegative principal symbol which is elliptic on the support of $\phi$. We construct $A$ as $A_1 + A_2$, where $A_2$ is constructed exactly as above, but $A_1$ is modified to reflect the hyperbolic rather than Euclidean structure at infinity. We shall take $A_1$ to be a zeroth order pseudodifferential operator in the 0-calculus of Mazzeo-Melrose \cite{MM}. Recall that the 0-calculus of pseudodifferential operators on a manifold with boundary is the natural class of pseudodifferential operators associated to differential operators generated by vector fields that \emph{vanish} at the boundary. This calculus is appropriate here since, if we compactify $C^2$ by adding circles at $r = \pm \infty$, with boundary defining functions $e^{\mp r}$, then the Laplacian $\Delta_{C^2}$ on $C^2$ is an elliptic combination of such vector fields. 

Using coordinates $(\rho, \omega)$ dual to $(r, \theta)$, we  define $A_1$ to be a zeroth order 0-pseudodifferential operator
with symbol $\zeta^2(r) (1 - e^{-\epsilon r}) \rho \cdot (1 + \sigma(\Delta_{C^2}))^{-1/2}$. In these coordinates, the symbol of $\Delta_{C^2}$ is 
$$
\sigma(\Delta_{C^2}) = \rho^2 + (\sech r)^2  \omega^2
$$
and the Hamilton vector field is 
$$
2 \Big( \rho \dbyd{r} + \tanh r (\sech r)^2 \omega^2 \dbyd{\rho} + \sech^2 r \ \omega \dbyd{\theta} \Big).
$$
Applying this to the symbol of $A_1$ gives the positive term
$$
2 \zeta(r)(1 + \sigma(\Delta_{C^2}))^{-1/2} \Big(2  \zeta'(r) (1 - e^{-\epsilon r}) \rho^2 +  \epsilon \zeta(r)  e^{-\epsilon r}\rho^2  + \tanh r (\sech r)^2 \zeta(r) (1 - e^{-\epsilon r}) \omega^2  \Big) 
$$
which is nonnegative everywhere and bounded below by a multiple of $\sigma(\Delta_{C^2})^{1/2}$ on the support of $\phi$. The rest of the proof is the same as in part (i), using the fact that zeroth order  operators in the 0-calculus are bounded on $L^2(C^2)$. 
\end{proof}

\begin{rem}\label{remark-conic} Exactly the same result holds if $X$ is replaced by an asymptotically conic manifold, with the same proof. We only have to replace $r^{-1} \ang{\xib}^{-1} \xb \cdot \xib$ in \eqref{a1} by the cotangent variable dual to $dr$. We shall use this remark in the next section.
\end{rem}

\begin{rem} We can rephrase this result as follows: the operators $T_X = \phi e^{-i t \Delta_X}$ and $T_{C^2}= \phi e^{-i t \Delta_{C^2}}$  , for $\phi \in C_c^\infty(M)$ supported away from the trapped set, and $T_{\tilde X} = \tilde \phi e^{-i t \Delta_{\tilde X}}$, for $\tilde \phi \in C_c^\infty(\tilde X)$ are bounded from $L^2(X)$ to $L^2([0,T]); H^{\demi}(X)$, resp. $L^2(C^2)$ to $L^2([0,T]); H^{\demi}(C^2)$, resp. $L^2(\tilde X)$ to $L^2([0,T]); H^{\demi}(\tilde X)$. 
\end{rem}

We return to the proof of Theorem~\ref{example}. Consider the function $u_e = (1 - \chi) u$. We can regard it as a function on $\tilde Z$, and as such it satisfies on the time interval $[0,T]$
\begin{multline}
(i \partial_t - \Delta_{\tilde Z}) u_e = w \equiv -2 \nabla \chi \cdot \nabla u + (\Delta_{\tilde Z} \chi) u  ; 
\qquad
u_e \big|_{t=0} = (1 - \chi) u_0 \in L^2(\tilde Z). 
\end{multline}
By Lemma~\ref{T}, $w \in L^2([0,T]; H^{-\demi}(\tilde Z))$. Let us write $u_e = u_e' + u_e''$, where $u_e'$ solves the PDE above with zero initial condition, and $u_e''$ solves the homogeneous equation $(i \partial_t - \Delta_{\tilde Z}) u_e'' = 0$ with initial condition $(1 - \chi) u_0$. By \cite{ST}, the Strichartz estimate \eqref{Strichartz-example} holds for $u_e''$. The function $u_e'$ is given by Duhamel's formula
$$
u_e'(\cdot, t) = \int_0^t e^{-i (t-s) \Delta_{\tilde Z}} w(\cdot, s) \, ds.
$$
We want to show that this is in $L^p_t L^q_x$ for Strichartz pairs $(p, q)$. Since we are dimension $d=2$, we have $p > 2$, and hence we can apply  the Christ-Kiselev Lemma \cite{CK}, which tells us that it is sufficient to show boundedness of the operator 
$$
w \mapsto \int_0^1 e^{-i (t-s) \Delta_{\tilde Z}} w(\cdot, s) \, ds
$$
from $L^2([0,T]; H^{-\demi}(\tilde Z))$ to $L^p_t L^q_x$. But, defining $T_{\tilde Z}$ as above,  this is $e^{i t \Delta_{\tilde Z}} T_{\tilde Z}^* w$ (for any $\tilde \phi$ equal to $1$ on the support of $\nabla \chi$). By Lemma~\ref{T} and duality,  $T_{\tilde Z}^*$ maps
$L^2([0,T]; H^{-\demi}(\tilde Z))$ to $L^2(\tilde Z)$, while by \cite{ST}, $e^{i t \Delta_{\tilde Z}}$ maps $L^2$ to $L^p_t L^q_x$. This shows that \eqref{Strichartz-example} holds for the function $u_e$. 

It remains to consider $u_i = \chi u$. We regard $u_i$ as a function on the hyperbolic cylinder $C^2$ since it is supported in the region where the metric is hyperbolic, and as such it satisfies on the time interval $[0,T]$
\begin{multline}
(i \partial_t - \Delta_{C^2}) u_i = w \equiv -2 \nabla \chi \cdot \nabla u + (\Delta_{C^2} \chi) u  ;
\qquad
u_i \big|_{t=0} = \chi u_0 \in L^2(C^2). 
\end{multline}
As we have seen, $w \in L^2([0,T]; H^{-\demi}(C^2))$. Let us write $u_i = u_i' + u_i''$, where $u_i'$ solves the PDE above with zero initial condition, and $u_i''$ solves the homogeneous equation $(i \partial_t - \Delta_{C^2}) u_i'' = 0$ with initial condition $ \chi u_0$. By Theorem~\ref{stric}, the Strichartz estimate \eqref{Strichartz-example} holds for $u_i''$. The function $u_i'$ is given by Duhamel's formula
$$
u_e'(\cdot, t) = \int_0^t e^{-i (t-s) \Delta_{C^2}} w(\cdot, s) \, ds.
$$
We want to show that this is in $L^p_t L^q_x$ for Strichartz pairs $(p, q)$. We use the Christ-Kiselev trick again and show that the operator
$$
w \mapsto \int_0^1 e^{-i (t-s) \Delta_{C^2}} w(\cdot, s) \, ds
$$
is bounded from $L^2([0,T]; H^{-\demi}(\tilde Z))$ to $L^p_t L^q_x$. But, defining $T_{C^2}$ as above,  this is $e^{i t \Delta_{C^2}} T_{C^2}^* w$ (for any $\tilde \phi$ equal to $1$ on the support of $\nabla \chi$). By (ii) of Lemma~\ref{T},  $T_{C^2}^*$ by duality maps
$L^2([0,T]; H^{-\demi}(C^2))$ to $L^2(C^2)$, while by Theorem~\ref{stric}, $e^{i t \Delta_{\tilde Z}}$ maps $L^2$ to $L^p_t L^q_x$. This shows that \eqref{Strichartz-example} holds for the function $u_i$, and completes the proof of Theorem~\ref{example}.

\end{proof}


\section{Asymptotically Euclidean (or conic) manifolds with filamentary hyperbolic trapped set}

In the previous section, we used the dispersive estimate from section \ref{hypman} for constant negative curvature manifolds to prove  
Strichartz estimates without loss. It is natural to ask if this result can be generalized to a variable negative curvature setting. 
In this section, we shall show that a more general class of manifolds 
with hyperbolic trapped set has this `Strichartz without loss' property.
The class of manifolds we will consider are asymptotically Euclidean (and more generally asymptotically conic) 
but the projection of their trapped set is contained in an open set where the metric has (variable) negative curvature, so that the  flow is hyperbolic there, and we will assume, as in \cite{NZ}, that the topological pressure $P(s)$ of the unstable Jacobian on the trapped set satisfies  $P(1/2)<0$ --- this last condition roughly means that the trapped set is thin enough, also called \emph{filamentary}, although it may contain an infinite number of closed geodesics.

An asymptotically conic manifold (or scattering manifold in the sense of \cite{Me}) 
is a complete non-compact Riemannian manifold $(M,g)$ which is the interior of a smooth compact manifold with boundary $\bbar{M}$ and such that a collar neighbourhood of the boundary is isometric to
\[ \Big( [0,\eps)_x\x \pl\bbar{M}, \frac{dx^2}{x^4}+\frac{h(x)}{x^2}\Big)\]
where $h(x)$ is a one-parameter family of metrics on $\pl\bbar{M}$ depending smoothly on $x\in [0,\eps)$.
Here $\pl\bbar{M}$ has really to be considered as the `points at infinity' of $(M,g)$.
The function $x$ can be extended to a nonnegative smooth function on $\bbar{M}$ and
the function $r=1/x$ is analogous to the radial function on Euclidean space. 

The geodesic flow $\Phi^t$, $t \in \RR$, is the flow of the Hamiltonian vector field $V_H$ associated to $H\in C^{\infty}(T^*M)$
defined by $H(m,\xi):=|\xi|^2_g$. The trapped set $K$ is defined by $K:=\Gamma^+\cap \Gamma^-$ where
\begin{equation}
\Gamma^{\pm}:=\{(m,\xi)\in T^*M \mid \Phi^t(m,\xi)\not\to\infty, t\to \mp\infty \}
\subset T^*M.
\label{Gamma}\end{equation}
Let us denote by $\pi: T^*M\to M$ the projection on the base, and let $d = \dim M$. 

The geodesic flow is said to be  hyperbolic on $U\subset S^*M$ if  
for all $m \in U$, the tangent space at $m$ splits into flow, unstable and stable subspaces such that 
\begin{equation} \label{splitting}
\begin{array}{rl}
i) & T_mS^*U=\rr V_H(m)\oplus E_m^+\oplus E_m^-, \quad \dim E_m^\pm=d-1\\
ii) & d\Phi_m^t(E_m^\pm)=E_m^\pm, \quad \forall t\in\rr\\
iii) & \exists \la>0, \quad ||d\Phi^t_m(v)||\leq Ce^{-\la|t|}||v||, \quad \forall v\in E_m^{\mp}, \ \pm t\geq 0.
\end{array}
\end{equation} 
for some uniform $\la>0$; here the norm can be taken with respect
to the Sasaki metric on the cotangent bundle (see \cite{Paternain}, Definition 1.17). This is true in particular for $U=S^*M$ if $M$ is a complete manifold with 
negative sectional curvatures contained in an interval $[-k_1,-k_0]$ for some $k_i>0$, see for instance
\cite[Th. 3.9.1]{Klin}.  

We define the unstable Jacobian $J^u_t(m)$ and the weak unstable Jacobian $J^{wu}_t(m)$ for the flow $\Phi^t$ at the point $m$ to be 
\begin{equation}\begin{aligned}
J^u_t(m) &= \det \big( d\Phi^{-t}(\Phi^t(m)) |_{E^+_{\Phi^t(m)}} \big), \\
J^{wu}_t(m) &= \det \big( d\Phi^{-t}(\Phi^t(m)) |_{E^+_{\Phi^t(m)} \oplus \RR V_H(m)} \big)
\end{aligned}\label{Ju}\end{equation}
where the volume form on $d$ dimensional subspaces of $T(T^*M)$ is induced by the Sasaki metric.
It follows from (iii) of \eqref{splitting} that $J^u_t(m), J^{wu}_t(m) \leq e^{-\lambda t}$ for $t > 0$.

If the geodesic flow on the trapped set $K$ is hyperbolic, and $s : K \to \RR$ is a continuous function, 
then  the \emph{topological pressure} of the unstable Jacobian at $s$ is a real number $P(s)$, whose definition is given by \eqref{pressuredef}. The topological pressure of the unstable Jacobian can be viewed as a real function $P$ of $s$. The quantity $P(0)$ is known as the topological entropy of the flow. 
For positive $s$, $P(s)$ in a sense measures two competing effects of the flow: the density of $K$ (the denser $K$, the longer points near $K$ stay close by under the flow) and the instability of the flow (the more unstable, the more quickly points near $K$  move away from $K$ under the flow). In our analysis we encounter products of square roots of the unstable Jacobian, which in view of \eqref{Z} and \eqref{pressuredef} make it natural to consider the topological pressure at $s=1/2$; if $P(1/2) < 0$ then the instability dominates, which is crucial in our main estimate (Lemma~\ref{mainest}) of this section. The first use, to our knowledge, of topological pressure in analytical estimates was by Nonnenmacher-Zworski \cite{NZ} following work of Gaspard-Rice \cite{GR} in the physics literature.

Our assumptions on $(M, g)$ and on the trapped set $K$ in this section are 
\begin{itemize}
\item[(A1)]  $(M,g)$ is asymptotically conic. 
\item[(A2)] There is an open set $X_- \subset X$ containing $\pi(K)$ which can be extended to a complete manifold $\tilde M$  with sectional curvatures bounded above by a negative constant (in particular, $M_-$ has  sectional curvatures bounded above by a negative constant). 
\item[(A3)] $M_-$ is geodesically convex in $M$: i.e. any geodesic entering $\pi^{-1}(M \setminus M_-)$ from $\pi^{-1}M_-$ remains in this region thereafter.
\item[(A4)] The topological pressure $P(s)$ of $K$ evaluated at $s=1/2$ is negative:
\begin{equation}
P\big(\frac1{2}\big) < 0.
\label{negpress}\end{equation}
\end{itemize}

For examples, see Section~\ref{examples}. 

\begin{rem} 
With these assumptions, the geodesic flow is hyperbolic on $S^*\tilde M$ and $K$
is the trapped set on both $S^*\til{M}$ and $S^*M$. On $S^*M$, the splitting satisfying \eqref{splitting}
only makes sense at points of $K$, but we can still consider the splitting 
$TS^*M_-=\rr V_H(m)\oplus E_m^+\oplus E_m^-$ coming from the inclusion $S^*M_- \subset S^*\til{M}$:
in particular, for all $m\in S^*M_-$ and all $t$ such that $\Phi^t(m)\in \pi^{-1}(M_-)$ we have 
\[\begin{gathered}
  d\Phi_m^t(E_m^\pm)=E_m^\pm,\\
 \exists \la>0, \quad ||d\Phi^t_m(v)||\leq Ce^{-\la|t|}||v||, \quad \forall v\in E_m^{\mp}, \textrm{ if }\pm t\geq 0.
 \end{gathered}\]   
 \end{rem}

\begin{rem} It seems likely that (A2) actually follows from (A1) and (A3); that is, that if $M_-$ is negatively curved and geodesically convex, then  we expect that it can always be extended to a complete manifold with negative sectional curvature. We do not pursue this question further here as it is a purely differential-geometric question.
\end{rem}
 
Our main result in this section is

\begin{theo}\label{sharpstrich} Let $(M, g)$ satisfy assumptions (A1) --- (A4) above. Then local-in-time Strichartz estimates without loss hold for $M$:  there exists $C>0$ such that
\begin{equation}
||e^{it\Delta}u_0||_{L^p((0,1),L^q(M))}\leq C||u_0||_{L^2(M)}
\label{sharpstrich3}\end{equation}
for all $u_0\in L^2(M)$ and $(p,q)$ satisfying \eqref{firstST}.
\end{theo}

\begin{rem} 
In dimension $d=2$, the condition $P(1/2)<0$ is equivalent to  the Hausdorff dimension of the trapped set satisfying  $d_{H}(K)<3$, or equivalently 
$d_H(K\cap S^{*}M)<2$. Note that this is the natural generalization  of the condition $\delta<n/2$ in our hyperbolic quotients examples above, since $d_{H}(K\cap S^*M)=2\delta+1$ in that case (recall $d=n+1$).
\end{rem}

\subsection{Topological pressure}\label{tp}
We now define the topological pressure of the flow $P(s)$ on the trapped set, following \cite{NZ} (which follows from Def. 20.2.1 of \cite{HK}): a set $E\subset K\cap S^*M$ is said to be $(\eps,T)$ separated if given $(x_1,v_1)\not=(x_2,v_2)$ in $E$, there exists $t\in[0,T]$ for which the distance between $\Phi^t(x_1,v_1)$ and $\Phi^t(x_2,v_2)$ is at least $\eps>0$.
For any $s \in \RR$ we define\footnote{For later convenience we define the topological pressure using the weak unstable Jacobian rather than the unstable Jacobian which is more standard. However, making this change only changes $\log Z_T(\epsilon, s)$ by $O(1)$, uniformly in $\epsilon$, and thus leads to the same value of $P(s)$.}
\begin{equation}
Z_T(\eps,s):=\sup_{E}\sum_{m\in E} (J^{wu}_T(m))^{s}
\label{Z}\end{equation}
where the $\sup$ is taken over all sets $E$ which are $(\eps,T)$-separated.
The pressure of $s$ is
\begin{equation}\label{pressuredef}
P(s)=\lim_{\eps\to 0}\limsup_{T\to\infty}\frac{1}{T}\log Z_T(\eps,s)
\end{equation} 

For instance, if the metric has constant curvature, one has
\begin{lem}\label{lempressure}
Let $(M, g)$ be a convex cocompact hyperbolic manifold of dimension $n+1$ with limit set of dimension $\delta$. Then the topological pressure at $s=1/2$ is given by $P(1/2) = \delta - n/2$.
\end{lem}
\begin{proof}
For a constant curvature $-1$ manifold, $J^u_t(m) = e^{-tn}$. It follows from   \eqref{Z} and \eqref{pressuredef} (and the footnote) that $P(s) = P(0) - ns$. 
But $P(0)$, which is the topological entropy, is equal to $\delta$  by a result of Sullivan \cite{Su}. 
\end{proof}

Finally, we recall the alternate definition 
of topological pressure given in \cite[Sec. 5.2]{NZ}, which turns out to be easier to use. If $\mc{V}=(V_{b})_{b\in B}$ is an open finite cover of $K\cap S^*M$, let $\mc{V}^T$ ($T\in\nn$) be the refined cover made of $T$-fold intersections
\[V_{\beta}:=\bigcap_{k=0}^{T-1}\Phi^{-k}(V_{b_k}), \quad \beta:=b_0b_1\dots b_{T-1}\in B^T,\] 
and consider the set $\mc{B}'_{T}\subset B^T$ of $\beta$ such that $V_\beta\cap K\not=\emptyset$.
For any $W\subset S^*M$ with $W\cap K\not=\emptyset$, define the coarse-grained 
unstable Jacobian 
\begin{equation}
S^K_{T}(W):=  \sup_{m\in W\cap K}\log J^u_T(m) = 
-\inf_{m\in W\cap K}\log \det(d\Phi^T(m)|_{E^+_m\oplus\rr V_H(m)}).
\label{stowa}\end{equation}
The topological pressure is defined by 
\[P(s):=\lim_{{\rm diam}\mc{V}\to 0}\lim_{T\to \infty}\frac{1}{T}\log 
\inf\Big\{ \sum_{\beta\in\mc{B}_T}\exp(sS^K_T(V_\beta)); \mc{B}_T\subset\mc{B}_T',
K\cap S^*M\subset\bigcup_{\beta\in\mc{B}_T}V_\beta\Big\}.\]
In particular, for all $\eps_0>0$ small, there exists $\eps_1>0$, such that for all $\eps<\eps_1$ and 
all covers $\mc{V}$ of $K\cap S^*M$ as above with diameter smaller than $\eps$, 
there is a $T_0\in\nn$, a set $\mc{B}_{T_0}\subset \mc{B}_{T_0}'$ such that $\{V_{\beta},\beta\in\mc{B}_{T_0}\}$
is an open cover of $K\cap S^*M$ and 
\begin{equation}
\sum_{\beta\in \mc{B}_{T_0}}\exp(sS^K_{T_0}(V_\beta))\leq \exp(T_0(P(s)+\eps_0/2)).
\label{T0}\end{equation}
Moreover $V_\beta$ are all included in $\pi^{-1}(M_-)$ since they are $\eps$ close to $K$.
Since by the chain rule one has 
$$
J^u_{T_0}(m) = \prod_{j=1}^{T_0} J^u_1(\Phi^{j-1}(m)),
$$
and since the unstable foliation is $\gamma$-H\"older \cite{SS}, we deduce that\footnote{Notice that $E^+_m$ makes sense since $\pi(V_\beta)$ 
is in the negative curved part $M_-$.} for any $V_{\beta}$ with $\beta\in \mc{B}_{T_0}$
\[|S^K_{T_0}(V_\beta)-S_{T_0}(V_\beta)|\leq \exp(CT_0\eps^\gamma), \textrm{ where }S_{T_0}(V_\beta):= \sup_{m \in V_\beta} \log J^u_{T_0}(m)
\]
for some constant $C$ depending only on $\eps_0$. 
Therefore, renaming the family $(V_{\beta})_{\beta\in\mc{B}_{T_0}}$ by
$(W_{a})_{a\in A_1}$, we get, by taking for instance $C\eps^{\gamma}\leq \eps_0/2$ 
\begin{equation}\label{suma1}
\sum_{a\in A_1}\exp(sS_{T_0}(W_a))\leq \exp(T_0(P(s)+\eps_0)).
\end{equation}

\subsection{Some examples}\label{examples}
We first give examples of Riemannian manifolds satisfying assumptions (A1) --- (A4).  
Consider any convex co-compact hyperbolic manifold $(M, g)$. Near infinity, it is conformally compact. That is, it admits a compactification to a compact manifold $\bbar{M}$ with a boundary defining function $x$, such that $x^2 g$ is a smooth nondegenerate metric up to the boundary of $\bbar{M}$. In other words, near infinity $g$ takes the form
$g = h/x^2$ for some smooth metric $h$ on $\bbar{M}$. 
We shall now modify the metric $g$ near infinity to a metric that is asymptotically conic, in such a way that the trapped set is left unchanged. (We recall that both asymptotically conic and conformally compact metrics are nontrapping near infinity.) This is straightforward: we have near infinity, in suitable coordinates $(x, y)$ where $x$ is a boundary defining function for $\bbar{M}$ and $y$ is a local coordinate on $Y = \partial \bbar{M}$
$$
g = \frac{dx^2 + h(x)}{x^2}.
$$
Here $h(x)$ is a smooth family of metrics on $Y$, i.e. is smooth in all its arguments. Changing variables to the `geometric' coordinate $r = \log(1/x)$, this reads
$$
g = dr^2 + e^{2r} h(e^{-r}).
$$
Assume that this is valid for $r \geq R$. Then, for some $R' \geq R$ we choose a function $f(r)$ such that $f(r) = e^{r}$ for $r \leq 2R'$ and $f(r) = c r$ for $r \geq 4R'$, and such that $f'(r) > f(r)/2r$ for all $r$. This is possible for all $R' \geq R$ and some $c$ depending on $R'$.
Define the metric 
\begin{equation}
g_{ac} = dr^2 + f(r)^2 h(e^{-r}), \ r \geq 2R; \quad g_{ac} \equiv g \text{ for } r \leq 2R.
\end{equation}
This is an asymptotically conic metric on $M$ and satisfies assumption (A2). The symbol of the Laplacian with respect to this metric is $\rho^2 + f(r)^{-2} (h^{-1}(e^{-r}))_{ij} \eta_i \eta_j$ and along geodesics we have 
$$
\ddot r = 2 \Big( \frac{f'(r)}{f(r)^3} |\eta|^2_{h^{-1}(e^{-r})} + \frac{e^{-r}}{f(r)^2} |\eta|^2_{\dot h^{-1}(e^{-r}) } \Big).
$$
Here $\dot h^{-1}(e^{-r})$ means $d/ds ( h^{-1}(s)) \mid s = e^{-r}$. 
The metric $d/ds (h^{-1}(s))$ is bounded above by a constant times $h^{-1}(s)$ uniformly for $s \in [0, \log R^{-1}]$. Also, using $f'(r) > f(r)/2r$, we see that $f'(r)/f(r)^3 \gg e^{-r}/f(r)^{2}$ for large $r$. It follows that for $R'$ sufficiently large, and $r \geq R'$, we have $\ddot r \geq 0$ and hence there is no trapped set in $r \geq R'$ for the metric $g_{ac}$; moreover,  the set $\{ r \leq R' \}$ is geodesically convex. Hence the metric satisfies condition (A3).\\ 

Finally to verify assumption (A4) it suffices to use Lemma \ref{lempressure}.

\subsection{Strategy of the proof}\label{strategy}
The proof of Theorem~\ref{sharpstrich} is much more involved than that of Theorem~\ref{example}. We will need to localize both in frequency and in time. To explain the idea, we first show how the Strichartz estimates on an asymptotically conic \emph{nontrapping} manifold $M$ may be deduced via frequency and time localization. Here we focus on  estimating the solution on a compact set contained in $M$. Thus, we consider $\chi u$, where  $\chi \in C_0^\infty(M)$ vanishes for small $x$. 

We introduce the semiclassical parameter $h$, where $h^{-1}$ will be (up to a constant) the frequency of our frequency-localized wave $u$. Let $\psi \in C_0^\infty(1/2, 2)$. Then the semiclassical dispersion estimate from \cite{BGT} says that for sufficiently small $c$,
\begin{equation}
\| \psi(h^2 \Delta_M) e^{-it\Delta_M} \|_{L^1(M) \to L^\infty(M)} \leq C t^{-n/2} \text{ for } t \in [0, ch].
\label{scde}\end{equation}
Let us assume that $u_0$ is localized near frequencies $\approx h^{-1}$ in the sense that $\psi(h^2 \Delta_M) u_0 = u_0$ (this will then be true for all times $t$). This assumption is harmless as a  Littlewood-Paley argument (see \cite{BGT}, Section 2.3.2 using the Littlewood-Paley estimate from \cite{Bo}) shows that if Strichartz is true for frequency-localized $u$, then it holds for all $u$. 
It follows from \eqref{scde} and Keel-Tao that Strichartz holds for $u$ on a time interval of length $ch$:
\begin{equation}
\| e^{-it\Delta_M} u_0 \|_{L^p[0,ch]; L^q(M)} \leq C \| u_0 \|_{L^2(M)}
\label{scs}\end{equation}
for $(p,q)$ satisfying \eqref{firstST}. Here, $c$ depends on the injectivity radius of $M$ and the support of $\psi$; for simplicity, below we assume that $c=1$. 

To extend this to a fixed length time interval, we use time cutoffs and local smoothing estimates. To define the time cutoffs, let $\varphi(s) \in C_0^\infty[-1,1]$ satisfy $\varphi(0) = 1$ and $\sum_{j \in \mathbb{Z}} \varphi(s - j) = 1$. Then we can write, for any Schr\"odinger wave $u(\cdot, t) = e^{-it\Delta_M} u_0$, $$\chi u = \sum_{j \in \mathbb{Z}}  \varphi(t/h - j)\chi  u \equiv \sum_j u_j,$$ where each $u_j$ is supported on a time-interval of length $2h$. We work on the time interval $[0,1]$ and assume that $h^{-1} = N \in \mathbb{N}$; thus, we need to consider $u_j$ for $j = 0, 1, \dots, N$. The functions 
$u_0$ and $u_N$ are dealt with from the semiclassical Strichartz estimate \eqref{scs}. So consider $u_j$ for $1 \leq j \leq N-1$. These functions satisfy the equation 
\begin{equation}
(i \partial_t - \Delta_M) u_j = h^{-1} \varphi'(t/h - j) \chi u + 2 \varphi(t/h - j) \big( \nabla \chi \cdot \nabla u - \Delta_M \chi u \big) \equiv w_j.
\label{wjdefn}\end{equation}
Since $M$ is nontrapping, we have from Lemma~\ref{T} and Remark~\ref{remark-conic} the local smoothing estimate 
\begin{equation}
\| \tilde \chi u \|_{L^2[0,1]; H^{1/2}(M)}  \sim h^{-1/2} \|  \tilde \chi u \|_{L^2[0,1]; L^2(M)}  \leq C \|  u_0 \|_{L^2(M)}
\label{lsnt}\end{equation}
for any $\tilde \chi \in C_0^\infty(M)$. Choose $\tilde \chi$ to be $1$ on the support of $\chi$. It follows that 
\begin{equation}
\sum_{j=1}^N \| w_j \|_{L^2_t; L^2(M)}^2 \leq C h^{-1} \|   u_0 \|_{L^2(M)}^2.
\label{wjsum}\end{equation}
We can express $u_j$ in terms of $w_j$ using Duhamel's formula:
\begin{equation}
u_j(t) =  \tilde\chi \int_{(j-1)h}^t e^{-i(t-s) \Delta_M} \tilde \chi w_j(s) \, ds.
\end{equation}
By the Christ-Kiselev lemma, if $p > 2$, in order to estimate the $L^p_t$ norm of $u_j$ in terms of the $L^2_t$ norm of $w_j$ it is sufficient to estimate the $L^p$ norm of $\tilde u_j$ defined by 
 \begin{equation}
\tilde u_j(t) = \int_{(j-1)h}^{(j+1)h} e^{-i(t-s) \Delta_M} w_j(s) \, ds = e^{-it\Delta_M} \int_{(j-1)h}^{(j+1)h} e^{is \Delta_M} w_j(s) \, ds.
\end{equation}
Now we can use the semiclassical Strichartz estimate since the time interval is $O(h)$. The dual estimate to \eqref{lsnt} gives 
 \begin{equation}
\| \int_{(j-1)h}^{(j+1)h} e^{is \Delta_M} w_j(s) \, ds \|_{L^2(M)}  \leq C h^{1/2} \| w_j \|_{L^2_t; L^2(M)}.
\label{aaa}\end{equation}
Then \eqref{scs} applied to this $L^2$ function shows that
\begin{equation}
\| \tilde u_j \|_{L^p_t; L^q(M)} \leq C h^{1/2} \| w_j \|_{L^2_t; L^2(M)}
\end{equation}
and the same estimate holds for $u_j$ by Christ-Kiselev. 
Squaring this inequality, summing over $j$ and using \eqref{wjsum} shows that 
\begin{equation}
\sum_{j=1}^{N-1} \| u_j \|_{L^p_t; L^q(M)}^2 \leq C \|  u_0 \|_{L^2(M)}^2.
\end{equation}
Using the continuous embedding from $l^2(\NN)$ to $l^p(\NN)$ if $p\geq 2$, we obtain 
\begin{equation}
\Big(\sum_{j=1}^{N-1} \| u_j \|_{L^p_t; L^q(M)}^p \Big)^{2/p} \leq C \|   u_0 \|_{L^2(M)}^2
\end{equation}
and this gives 
\begin{equation}
\|   u \|_{L^p[0,1]; L^q(M)} \leq C \|  u_0 \|_{L^2(M)}
\end{equation}
with $C$ independent of $h$. 

\

Now suppose that $M$ is trapping, but obeys assumptions (A1) -- (A4) above.  In that case, the local smoothing estimate definitely fails \cite{Burq1}, \cite{Doi}, and then the argument only gives Strichartz estimate with a loss (i.e. with additional negative powers of $h$ on the right hand side) arising from the loss in the local smoothing estimate. If the trapped set has negative topological pressure, the local smoothing loss is  $|\log h|^{1/2}$ as follows from work of Nonnenmacher-Zworski \cite{NZ} and Datchev \cite{Da}, Theorem~\ref{lsll} below (see also~\cite{Burq1, BZ}). 
In combination with the argument above, this gives Strichartz with logarithmic loss, as shown for example in the case of the exterior of several convex obstacles by Burq \cite{Burq1}. 

However, when the topological pressure is negative then more is true: essentially from the work of Anantharaman \cite{An} and Nonnenmacher-Zworski \cite{NZ}, it follows that the semiclassical Strichartz estimate can be \emph{improved} by a logarithm: it is valid not just on a time interval of length $O(h)$, but actually on an interval of length $h |\log h|$ --- see Theorem~\ref{scsl} below. Then it turns out that this logarithmic improvement exactly compensates for the logarithmic loss in the local smoothing estimate and we recover the Strichartz estimate on a fixed finite time interval without loss. This is achieved by localizing in time on intervals of length $h |\log h|$ rather than $h$ for the part of $u$ localized near the trapped set. 

More precisely, we proceed as  above but with $\varphi(t/h - j)$ replaced with $\varphi(t/h|\log h| - j)$; that is, we localize to time intervals of length $h |\log h|$ which is the maximum for which we can apply the semiclassical Strichartz estimate. 
We then write $w_j = w_j'+ w_j''$, where
\begin{equation}\label{decomegaj}
w_j' = \frac1{h|\log h|} \varphi'(t/h|\log h| - j) u_j, \quad w_j'' =   2 \varphi(t/h|\log h| - j) \big( \nabla \chi \cdot \nabla u - \Delta_M \chi u \big) .
\end{equation}
Then $w_j''$ is of size $\sim h^{-1}$ (since a derivative of $u$ costs $h^{-1}$). On the other hand, it is supported in the nontrapping region and may be dealt with as above, as the local smoothing estimate is valid without loss in the nontrapping region. 
The other term, $w_j'$, is supported in the trapped region but is of size $O((h|\log h|)^{-1})$. On this term we apply the local smoothing estimate, losing $|\log h|^{1/2}$ as compared to the argument above. When we apply the dual estimate at the step \eqref{aaa} we lose a further $|\log h|^{1/2}$, and then applying semiclassical Strichartz completes the argument with no overall loss. 
The details are given in  Section~\ref{sec-proof}. 

In summary, the key ingredients of the proof will be the following three results:

\begin{theo}[Strichartz in the nontrapping region]\label{htwt} Let $M$ be an asymptotically conic manifold, and suppose that $\chi \in C^\infty(M)\cap L^\infty(M)$ vanishes in a neighbourhood of $\pi(K)$ where $K \subset T^*M$ is the trapped set.  Then we have Strichartz estimates without loss for $\chi e^{-it\Delta_M} u_0$:
\begin{equation}
\| \chi e^{-it\Delta_M} u_0\|_{L^p[0,1]; L^q(M)} \leq C \| u_0 \|_{L^2(M)}
\label{gs}\end{equation}
for all $(p,q)$ satisfying \eqref{firstST}, $p > 2$.
\end{theo}

\begin{theo}[Local smoothing with logarithmic loss]\label{lsll} Suppose that $M$ satisfies assumptions (A1) --- (A4). Then for any $\chi \in C_0^\infty(M)$ and $\psi \in C_c^\infty(1/2, 2)$, we have 
\begin{equation}
\| \chi e^{-it\Delta_M} \psi(h^2 \Delta) u_0 \|_{L^2[0,1]; L^2(M)} \leq C (h|\log h|)^{1/2} \| u_0 \|_{L^2(M)}.
\label{smoothing}\end{equation}
Moreover, if $\chi$ is supported outside the trapping region, then the estimate holds without the logarithmic loss in $h$ on the right hand side. 
\end{theo}

\begin{theo}[Semiclassical Strichartz on a logarithmic interval]\label{scsl} Suppose that $M$ satisfies assumptions (A1) --- (A4). Then for any $\chi$ supported in $M_-$, we have on a time interval of length $h |\log h|$
$$
\| \chi e^{-it\Delta_M} \psi(h^2 \Delta) u_0 \|_{L^p[0,h |\log h|]; L^q(M)} \leq C \| u_0 \|_{L^2(M)}
$$
for all $(p,q)$ satisfying \eqref{firstST}.
\end{theo}

\subsection{Strichartz in the nontrapping region}
In this section we sketch how to prove Theorem~\ref{htwt}. The argument is quite related to similar ideas in~\cite{ST}. Here, we  follow fairly straightforward modifications of the argument in \cite{HTW} for nontrapping metrics. Notice that none of the results about `Local Schr\"odinger integral operators' in section 3 of \cite{HTW} use the nontrapping property, which only enters when the local smoothing estimate is used. To adapt the results of \cite{HTW} to prove \eqref{gs}, we modify the definition of the Banach space $X$ in (4.6) of that paper to include a cutoff function $\chi_+$, supported in the nontrapping region and equal to $1$ on the support of $\chi$ in \eqref{gs} in the $\| f \|_{H^{1/2, -1/2 - \rho/2}}(M)$ term. This cutoff function then needs to be included in Lemma 4.4 and  Lemma 5.5 of \cite{HTW}.  The proof in section 6 then goes through provided that $\psi_\alpha^{(0)}$ is supported in the nontrapping region.

\subsection{Local smoothing effect}\label{lse} 

Theorem~\ref{lsll} follows fairly directly from the resolvent estimate from \cite{Da} (which generalizes \cite[Theorem 5]{NZ} to scattering manifolds):
\begin{equation}\label{eq.estsemi}\| \chi ( h^2\Delta -(1\pm i\epsilon))^{-1} \chi \|_{L^2 \rightarrow L^2}\leq C \frac {|\log h|} { h}, \qquad 0 < h < h_0 < \! < 1,
\end{equation}
with $C$ independent of $\epsilon$ (and $h$). 
From this, we deduce~\eqref{smoothing}, following~\cite{BGT2}. Indeed, denote by $T$ the operator 
$$ T= \chi e^{it\Delta}\psi(h^2\Delta).$$
The boundedness of $T$ from $L^2$ to $L^2((-\infty,+\infty);L^2(M))$ is equivalent to the boundedness of its adjoint $T^*$ from $L^2((-\infty,+\infty);L^2(M))$ to $L^2$ (with same norm), which in turn is equivalent to the boundedness of $TT^*$ from  $L^2((-\infty,+\infty);L^2(M))$ to itself (with same norm squared). But
$$ TT^* f(t)= \int_{- \infty} ^{+ \infty}  \chi e^{i(t-s) \Delta}\psi^2(h^2\Delta)\chi  f(s) ds= \chi\int_{- \infty}^t+ \chi \int_t^{+\infty} \equiv \chi A_1 f(t) + \chi A_2 f(t)$$
and it is enough to estimate for example $\chi A_1f $. 
For this we can assume that $f$ has compact support, and consider $u_\epsilon= e^{-\epsilon t }A_1f(t)$ and $f_\epsilon =e^{-\epsilon t }f $ (notice that $u_\epsilon$ is supported in the set $\{t \geq C\}$) which satisfy
$$ (i \partial _t + \Delta +i \epsilon ) u_\epsilon = \psi^2(h^2\Delta)\chi f_\epsilon$$
Taking Fourier transform with respect to the variable $t$, we get
$$ \chi \widehat{u_{\epsilon}} (\tau)= \chi ( \Delta -(\tau -i \epsilon))^{-1} \psi^2(h^2\Delta)\chi \widehat{ f_\epsilon}$$
and according to the Plancherel formula (recall that the Plancherel formula is true for functions taking values in any (separable) Hilbert space), and using~\eqref{eq.estsemi} we obtain
$$ \|\chi u_\epsilon \|_{L^2((-\infty,+\infty);L^2(M))} \leq C h \log(1/h)\|\chi f_\epsilon \|_{L^2((-\infty,+\infty);L^2(M))}.$$
Letting $\epsilon>0$ tend to $0$ we obtain that the contribution of $\chi A_1$ to $TT^*$ satisfies the required estimate. The other contribution $\chi A_2$ is dealt with similarly.

The improved estimate when the support of $\chi$ does not meet the trapped set is a consequence of Lemma~\ref{T} and Remark~\ref{remark-conic}.

\subsection{Semiclassical Strichartz on a logarithmic interval}
Using the work of Nonnenmacher-Zworski \cite{NZ} (which follows techniques of Anantharaman \cite{An}), we shall obtain a sharp dispersive estimate for the propagator $e^{-it\Delta_M}$ for time $t\in (0,h|\log h|)$ and in frequency localization in windows of size $h$.  
We want to prove the following:
\begin{prop}\label{displogh}
There exists $\delta>0$ and $C>0$ such that for all $\psi\in C_0^\infty((1-\delta/2,1+\delta/2))$, all $t\in (0,h|\log h|)$ with $h\in(0,h_0)$ small, we have for every $\chi \in C_0^\infty(M)$ supported in $M_-$
\[||\chi e^{it\Delta}\psi(h^2\Delta)\chi||_{L^1\to L^\infty}\leq Ct^{-n/2}.\]
\end{prop}
Then Theorem~\ref{scsl} follows immediately from this by applying the main result of Keel-Tao \cite{KT}. 

The proof of Proposition~\ref{displogh} decomposes in several parts. Let us first introduce the objects taken from \cite{NZ} that we need to use for the proof.\\
   
Without loss of generality and to simplify notation, we assume as in \cite{An} that 
the injectivity radius of $M$ is larger than $1$.
For $t\in(0, h)$, the result is essentially contained in 
 \cite{BGT,BT} since we are localized in a compact set of $M$. 
Now take $s_0\in [0,1]$ and an integer $L$, with $1 \leq L\leq \log(1/h)$. We want to obtain a dispersive estimate for $U(L+s_0)$, where 
\[U(t):=e^{ith\Delta}, \quad h\in (0,h_0)\]
by following \cite{An,NZ}. 
We consider, as in Section 6.3 of \cite{NZ}, a microlocal partition of unity $(\Pi_a)_{a\in A}$ of 
the energy layer $\mc{E}^{\delta}:=\{(m,\xi)\in T^*M, |\xi|\in (1-\delta,1+\delta)\}$ for some $\delta>0$ small. 
Let us recall how the partition $(\Pi_a)$ is defined.
The operators $\Pi_a$ are associated to an open covering 
$(W_a)_{a\in A}$ of $\mc{E}^\delta$ in the sense that the semi-classical wavefront set $\wf(\Pi_a)\subset W_a$ and $\sum_{a\in A}\Pi_a=I$ microlocally near $\mc{E}^{\delta/2}$, i.e. 
$\wf(\Pi_{\infty})\cap \mc{E}^{\delta/2}=\emptyset$ if $\Pi_{\infty}$ is defined by $\Pi_{\infty}:=I-\sum_{a\in A}\Pi_a$. 
Following Section 5.2 and 5.3 in \cite{NZ}, the set $A$ is decomposed in $3$ parts, $A=A_1\sqcup A_2\sqcup \{0\}$. The open set $W_0$ is defined by 
\[W_0:=\mc{E}^\delta\cap \pi^{-1}( M \setminus M_-).\] 
$(W_a)_{a\in A_1}$ is chosen so that $W_a\subset M_-$ if $a\in A_1$ and, as $a$ ranges over $A_1$, these sets covers $K\cap\mc{E}^\delta$ in such a way that for any $\eps_0>0$ fixed small,
there exist $\delta,\eps>0$, $T_0\in\nn$ with $W_a\subset \{m\in T^*M,d(W_a,K)\leq \eps\}$ and 
\begin{equation}\label{epsilon}
\sum_{a\in A_1}\exp\Big(s S_{T_0}(W_a)\Big)\leq \exp\Big(T_0(P(s)+\eps_0)\Big),
\end{equation}
where $S_{T_0}(W_a)$ is defined by \eqref{stowa}. This is possible as explained in Section~\ref{tp} and using
the homogeneity of the Hamiltonian on $T^*M$ to deal with $\mc{E}^\delta$ instead of $S^*M$. Finally, the $W_a$ for  $a \in A_2$ are defined so that there exists $d_1>0$ such that 
\[d(W_a,\Gamma^+\cap\mc{E}^\delta)+d(W_a,\Gamma^-\cap \mc{E}^\delta)>d_1\] 
where $\Gamma^{\pm}$ are the forward/backward trapped sets defined in \eqref{Gamma}. 
By \cite[Lem. 5.1]{NZ}, there exists $L_0\in\nn$ such that for all $a\in A_2$
\begin{equation}
\Phi^{t}(W_a)\subset W_0 \textrm{ for }t\geq L_0 \textrm{ or }t\leq -L_0.
\end{equation}
For notational simplicity we replace both $T_0$ and $L_0$ by $\max(T_0, L_0)$; hence we have 
\begin{equation}\label{defN0}
\Phi^{t}(W_a)\subset W_0 \textrm{ for }t\geq T_0 \textrm{ or }t\leq -T_0.
\end{equation}

To summarize, the energy layer $\mc{E}^{\delta}$ is decomposed into the part $W_0$ covering the (spatial) 
infinity of $\mc{E}^\delta$,
the part $\cup_{a\in A_1}W_a$ covering the trapped set $K\cap \mc{E}^{\delta}$, and finally the part  
covering the complementary, whose flowout by $\Phi^t$ lies in $W_0$ after some large (positive or negative) time.\\

We write $U(t)=e^{ith\Delta}$ and we shall prove 
\begin{equation}
||\chi \psi(h^2\Delta)U(T)u||_{L^\infty}\leq C(Th)^{-n/2}||u||_{L^1}
\label{wts}\end{equation}
for $T\in (2T_0,|\log h|)$. 
First we need a technical lemma.

\begin{lem}\label{tech} (i) Let $(M, g)$ be a scattering  manifold and $\Delta_M$ its Laplacian. Then for any $\psi \in C_0^\infty(\RR)$, $\psi(h^2 \Delta_M)$ is a semiclassical scattering operator of order $(-\infty, 0, 0)$ in the sense of Wunsch-Zworski \cite{WZ}. 

(ii) Let $\psi$ be as above, then for each fixed $t$, the operator $\psi(h^2 \Delta_M) e^{-ith\Delta_M}$ is a semiclassical Fourier Integral operator  associated to the canonical relation
$$
\{ \big((z, \zeta), (z', \zeta')\big) \in T^*M \times T^*M \mid (z, t\zeta) = \exp (z', t\zeta') \}
$$
\end{lem}

\begin{proof}  
(i) This follows from the argument in \cite{HV} for scattering pseudodifferential operators. (We also remark that a similar, weaker result proved under more general assumptions about the nature of the ends of the manifold by Bouclet \cite{Bo} would also suffice for our purposes.)

(ii) It is shown in \cite{NZ} that $e^{ith\Delta_h \phi(h^2 \Delta)}$ is a semiclassical FIO for each fixed $t$ and all $\phi \in C_0^\infty(\RR)$. Thus, using the result of (i), $\psi(h^2 \Delta) e^{ith\Delta \phi(h^2 \Delta)}$ is a semiclassical FIO. If $\phi = 1$ on the support of $\psi$ then this is precisely $\psi(h^2 \Delta) e^{ith\Delta_h}$ by functional calculus, proving the result.
\end{proof}

We decompose $T > 2T_0$ in the form $T = L-1 + s_0 = (2+N)T_0 + t_0$, where 
$L, N \in \nn$ and $s \in (0, 1]$,  $t_0 \in (0, T_0]$. Choosing $\psi_+ \in C_c^\infty(\RR)$ to be $1$ on the support of $\psi$, we have 
$$
\psi(h^2 \Delta) U(T) = \psi(h^2 \Delta) \psi_+^{L} (h^2 \Delta) U(T) =    \psi(h^2 \Delta) U(s_0) \Big( \psi_+(h^2 \Delta) U(1) \Big)^{L-1} \psi_+(h^2 \Delta)
$$
and we can decompose  
\begin{equation}
\psi_+(h^2 \Delta) U(1)=\sum_{a\in A\cup \infty}U_a, \quad U_a:=\psi_+(h^2 \Delta) U(1)\Pi_a.
\label{Ua}\end{equation}
Hence we may write 
\begin{equation}
\psi(h^2 \Delta) U(T)= \sum_{\alpha\in A^L}\psi(h^2 \Delta) U(s_0)\Pi_{\alpha_L}U_{\alpha_{L-1}} \dots U_{\alpha_1}\psi_+(h^2 \Delta) +R_{T}
\label{UT-decomp}\end{equation}
where $R_T$ term is the sum over all sequences $\alpha$ containing at least one index $\alpha_j=\infty$.  
The first estimate one obtains corresponds to Lemma 6.5 in \cite{NZ}: 
\begin{lem}\label{lemrn} 
If $\psi,\psi_+(x,hD)$ are chosen as above,  one has  
\[ \| \chi R_T \chi \|_{L^1\to L^\infty} =O(h^{\infty})\]
for $1 \leq T\leq |\log h|$, with implied constants independent of $T$.
\end{lem}

\begin{proof} Since both $\Pi_\infty$  and  $\psi_+(h^2 \Delta)$ are order zero pseudodifferential operators and they have disjoint operator wavefront set,  the composition $\Pi_\infty \psi_+(h^2 \Delta)$ is $O(h^\infty)$ as an operator from $L^2$ to $L^2$. The other factors in \eqref{UT-decomp} are all bounded from $L^2$ to $L^2$, and there are at most $C e^{CT}$ terms in the sum. As $T \leq |\log h|$ this contributes at most a factor of a fixed power of $h$. Hence $\| R_T \|_{L^2 \to L^2} = O(h^\infty)$.

To get an $L^1 \to L^\infty$ estimate from this, we compose on the left with $\psi(h^2 \Delta_M)$ and observe that we still obtain an $O(h^\infty)$ estimate if we pre- and post-multiply by $(1 + \Delta_M)^m$ for any $m$, since this has the effect of increasing the operator norm by at most $C h^{-4m}$. This is equivalent to an $O(h^\infty)$ estimate from Sobolev spaces $H^{-m}(M)$ to $H^m(M)$, from which we obtain $L^1 \to L^\infty$ by Sobolev embedding for $m$ larger than half the dimension of $M$. 
\end{proof}

The second estimate needed is similar to Lemma 6.6 of \cite{NZ}, but it is even better since 
we cut on the left on a compact set.  Let $\mc{A}_L\subset (A\setminus\{0\})^{L}$ defined by
\begin{equation}\label{defAN}
\alpha\in \mc{A}_L\iff\left\{
\begin{array}{l}
\Phi^{1}(W_{\alpha_j})\cap W_{\alpha_{j+1}}\not=\emptyset, \quad j=1,\dots,L-1\\
\textrm{ and }\alpha_j\in A_1 \textrm{ for all } j=T_0,\dots,L-T_0
\end{array}
\right. . \end{equation}

\begin{lem}\label{lemAN}
If $\alpha\in A^{L}\setminus \mc{A}_L$ and $\chi,\psi,\psi_+(x,hD)$ are chosen as above, then 
\begin{equation}\label{a2}
||\chi\psi(h^2\Delta) U(s_0)\Pi_{\alpha_L}U_{\alpha_{L-1}}\dots U_{\alpha_1}\psi(x,hD)\chi||_{L^1\to L^\infty}=O(h^{\infty})
\end{equation}
for $2T_0 < T\leq |\log h|$, with $T=L-1+s_0$ and $s_0\in[0,1]$.
\end{lem}

\begin{proof} This is proved in the same way as the previous lemma. We only need to show that each term in \eqref{UT-decomp} corresponding to a multi-index $\alpha \in A^{L}\setminus \mc{A}_L$ has a factor which is $O(h^\infty)$ as a map from $L^2$ to $L^2$. This follows directly from Egorov's theorem if there is a $j$ such that $\Phi^{1}(W_{\alpha_j})\cap W_{\alpha_{j+1}}=\emptyset$. Indeed, referring back to Lemma~\ref{tech} we can write 
$$U_{\alpha_{j+1}}\circ U_{\alpha_j} = \psi_+(h^2 \Delta) U(1) \Pi_{\alpha_{j+1}} \psi_+(h^2 \Delta) e^{ih\Delta \psi_{++}(h^2 \Delta)} \Pi_{\alpha_j}
$$
where $\psi_{++} \in C_c^\infty(\RR)$ is $1$ on the support of $\psi_+$. 
By Egorov, we have $e^{ih\Delta \psi_{++}(h^2 \Delta)} \Pi_{\alpha_j} = 
Q e^{ih\Delta \psi_{++}(h^2 \Delta)}$ for some pseudodifferential operator $Q$ with wavefront set given by $\Phi^{-1}(WF'(\Pi_{\alpha_j}))$. Since this is disjoint from the operator wavefront set of $\Pi_{\alpha_{j+1}}$ by hypothesis, this factor is $O(h^\infty)$ as a map from $L^2$ to $L^2$.

If either $\alpha_1 = 0$ or $\alpha_L = 0$ then the $O(h^\infty)$ estimate is immediate because $\Pi_0$ is microsupported in $M \setminus M_-$ and $\chi$ is supported in $M_-$. If any of the other $\alpha_j = 0$ then the $O(h^\infty)$ estimate follows because of assumption (A3), which implies that either $a_1 = 0$ or $a_L = 0$, or else the condition $\Phi^{1}(W_{\alpha_j})\cap W_{\alpha_{j+1}}=\emptyset$ has to hold for some intermediate $j$, showing that we are back in the situation considered above. Similarly, if $\alpha_j\in A_2$ for some $T_0 \leq j \leq L-T_0$ then \eqref{defN0} shows that we are again back in the situation considered above. 
\end{proof}

This Lemma clearly implies the bound 
\begin{equation}\label{sumA2}
\sum_{\alpha\notin\mc{A}_L}||\chi \psi(h^2\Delta)U(s_0)\Pi_{\alpha_L}U_{\alpha_{L-1}}\dots U_{\alpha_1}\chi||_{L^1\to L^\infty}=O(h^{\infty})
\end{equation}
since $|A|^{L}=O(h^{-\log |A|})$.

It remains to deal with the elements $\alpha\in\mc{A}_L$. We can obtain, again essentially from the analysis
of \cite{An} (and in a comparable way to \cite[Prop 6.3]{NZ}), the following bounds:

\begin{lem}\label{mainest}
Let $\chi$ be as above and let $\eps_0$ be
the small parameter in \eqref{epsilon}. Then for all small $\eps>0$, there exists $C>0$ such that  
\begin{equation}
\begin{gathered}
\sum_{\alpha\in \mc{A}_L}||\chi \psi(h^2\Delta)U(s_0)\Pi_{\alpha_L}U_{\alpha_{L-1}}\dots 
U_{\alpha_1}\chi||_{L^1\to L^\infty}\leq 
Ch^{-d/2}e^{T(P(1/2)+\eps_0+\eps)}
\end{gathered}
\label{aaaa}\end{equation}
for all $h\in (0,h_0)$ and $4T_0 \leq T\leq |\log h|$, where $T=L-1+s_0$ with $s_0\in[0,1]$.
\end{lem}

\textsl{Proof}:  We start by proceeding as in \cite[Sec. 3]{An}. 
If the cover is taken thin enough, we may use coordinates $(z,\xi)$ in each $W_a$, $a\in A_1$,
where $z\in\pi(W_a)$ and $\xi\in T^*_zM$ are cotangent variables. 
We can write for $u\in L^1(M)$ and $z\in \pi(W_{\alpha_1})$ 
\[\Pi_{\alpha_1} \chi u(z)=\int_{\pi(W_{\alpha_1})} \delta_{y}(z) u(y) \, dy +O(h^{\infty}), \textrm{ with }\]
\[\delta_{y}(z):= \frac{1}{(2\pi h)^d}\int_{(z,\xi)\in W_{\alpha_1}} e^{i\frac{(z-y)\xi}{h}}\sigma(z,\xi)d\xi.  
\]
where $\sigma(x,\xi)$ is the local symbol of $\Pi_{\alpha_1}\chi$ in $W_{\alpha_1}$. An upper bound for the left hand side of \eqref{aaaa} is then the sum over all $\alpha \in \mc{A}_{L}$ of 
\begin{equation}
\sup_{y, z}  \bigg| \Big( \big( \psi(h^2 \Delta) U(t)\Pi_{\alpha_{J+1}} \big) U_{\alpha_{J}} \dots U_{\alpha_2}e^{ih\Delta_0}\delta_{y} \Big)(z) \bigg|.
\label{bbb}\end{equation}
where we shall choose $t=s_0$ and $J=L-1$. Thus we take $\Pi_{\alpha_1} \chi u$ and evolve it through $e^{ih\Delta_0}$ then microlocally cutoff in $W_{\alpha_2}$, evolve again, microlocally cut off again, and so on. 
For Anosov flows, it is shown in \cite[Sec. 3]{An} that, for any $J,K\in\nn$ fixed (independently of $L$), there exists a function 
$S_J(.,t)\in C^{\infty}(\pi(W_{\alpha_{J}}))$ and $b_J(.,h,t)\in C^{\infty}(\pi(W_{\alpha_{J}}))$ 
with $b_J$ smooth in $h\in[0,h_0)$ such that for $t\in[0,1]$
\begin{equation} \begin{gathered} 
\Big( \big( \psi(h^2 \Delta) U(t)\Pi_{\alpha_{J+1}} \big) U_{\alpha_{J}} \dots U_{\alpha_2}e^{ih\Delta_0}\delta_{y} \Big) (z)=(2\pi h)^{-d/2}e^{\frac{iS_{J}(z,t)}{h}}b_{J,K}(z,h,t)+R_{J,K}(h,t) \\
b_{J,K}(z,h,t)=\sum_{k=0}^Kh^kb_{J;k}(z,t)
\end{gathered}\label{bJ;k}\end{equation}
with $||R_{J,K}(h,t)||_{L^2}=C_KJh^{K}$ for some $C_K>0$ uniform in $t,J$ (this estimate is shown in \cite[Lemma 3.2.2]{An}). 
The function $S_{{J}}(z,t)$ generates a smooth Lagrangian
submanifold $\mc{L}_{J+t} = \mc{L}_{J}(t)=\{(z,d_zS_{J}(z,t))\in T^{*}M; z\in \pi(W_{\alpha_{J}})\}$ which is part of the graph of the canonical transformation $\Phi^{J+t}$, namely that part with first coordinate lying in the Lagrangian $\{ (y, \xi) \mid 1 - \epsilon < |\xi| < 1 + \epsilon \}$. 

\begin{rem}\label{weak-unstable}
The key to the proof of Lemma~\ref{mainest} which we owe to \cite{An} is the following fact: as  $J \to \infty$ and since 
we only consider $\alpha \in \mc{A}_L$, the geodesics generating $\mc{L}_{J}(t)$ lie entirely within $\pi^{-1}(M_-)$ which has sectional curvatures bounded above by a negative constant, these Lagrangians $\mc{L}_{J+t}$ converge uniformly (indeed, exponentially) to the weak unstable foliation as $J \to \infty$. This will allow us to compare the size of $b_J$ with the weak unstable Jacobians $J^{wu}_t(m)$, as we shortly show. 
\end{rem}

We now take $J=L-1$, $t=s_0$ and $K$ large. We obviously have  $||\chi R_{L,K}(h,s_0)||_{L^2}\leq C_Kh^{K}|\log h|$
so using Sobolev embedding arguments as in the proof of Lemma \ref{lemrn}, 
we can replace the $L^2$ norm by the $L^\infty$ norm,
up to a loss of $h^{-d/2-1}$, so  
\[||\chi R_{L,K}(h,s_0)||_{L^\infty}\leq Ch^{K-d/2-1}|\log h|.\]
Taking $K$ large enough and summing over the $\alpha\in \mc{A}_L$, the number of which is bounded by $h^{-\log |A|}$, 
we conclude that these terms do not contribute.\\
   
It thus remains to study the $L^\infty$ norm of  elements of the form
\[(2\pi h)^{-d/2}\chi e^{\frac{iS_{L-1}(z,s_0)}{h}}b_{L-1,K}(z,h,s_0);\]
that is, the $L^\infty$ norm of $b_{L-1,K}$. We can essentially use the estimates in \cite{An} but first we need to make some remarks on the different partitions of unity used here as compared to \cite{An}. There, the quantum partition of unity is implemented by multiplication operators that cut off at the scale $h^{-\kappa}$, $0 < \kappa < 1/2$, while here we use semiclassical pseudodifferential operators with symbols smooth in $h$. There are two main differences: in \cite{An}, the multiplication operators are trivially bounded on $L^2(M)$ with operator norm $1$, while in our case, the operator norm of our microlocal cutoffs is $1 + O(h)$ since the principal symbols are bounded by $1$. This is inessential since it contributes at most a factor $(1 + Ch)^{|\log h|}$ to each term, which is bounded uniformly as $h \to 0$. Second, we need to replace the estimate on the derivatives of the microlocal cutoffs from $ |D^m A_a | \leq C h^{-m\kappa}$ in \cite{An} to $\| \ad_m(D, \Pi_a) \|_{L^2 \to L^2} \leq C$, where $D$ indicates differentiation and $\ad_m$ indicates the $m$th iterated commutator (which is even better than in \cite{An}, as we do not get any negative powers of $h$ in our case). 

With these remarks made, we can follow the analysis of Section 3.2  of \cite{An}. Let us define $J^{t}_{\mc{L}_L+s_0}(z)$ to be the Jacobian of the map $\Phi^{t}$, restricted to $\mc{L}_{L+s_0}$, and evaluated at
$\zbar = (z,dS_{L}(s_0)(z))\in\mc{L}_{L +s_0}$. Then the construction of \cite{An} shows that $b_{L-1;k}(z, s_0)$ is only nonzero if $\Phi^{-j}(z,dS_{L-1}(z,0))\in\pi(W_{\alpha_{L-j}})$
for all $j=1,\dots,L$ and $k=0,\dots, K$, in which case 
\begin{equation}
|b_{L-1;k}(z,s_0)|\leq C_kL^{3k}(J_{\mc{L}_{-(L-1+s_0)}}^{L-1+s_0}(\zbar))^{1/2}. 
\label{ccc}\end{equation}
Notice that this is the analogue of \cite[Lemma 3.2.1]{An}.
Let us write $T=T_1 + N T_0 + T_0$, where $T_1 = T_0 + t_0 \in [T_0, 2T_0]$. Then we can decompose 
\begin{equation}\label{Jterms}
\begin{gathered}
J_{\mc{L}_{L-1 +s_0}}^{-T}(\zbar)=J_{\mc{L}_{L-1+s_0}}^{-T_1}(\zbar) \times
J^{-T_0}_{\mc{L}_{(N+1)T_0}}\big(\Phi^{-T_1}(\zbar)\big)  \times
J^{-T_0}_{\mc{L}_{NT_0}}\big(\Phi^{-T_1-T_0}(\zbar)\big) \dots
\\ \quad \quad \dots\x  J^{-T_0}_{\mc{L}_{T_0}}\big(\Phi^{-T+T_0}(\zbar)\big).
\end{gathered}\end{equation}
The first and last Jacobian factors are uniformly bounded with respect to $L$; they only depend
on $T_0$ since they can be written as a supremum of the Jacobian of the flow at some time 
bounded by $2T_0$ on some set independent of $L$. Now using Remark~\ref{weak-unstable},  
by assuming that $T_0$ is large enough, the Lagrangians $\mc{L}_{jT_0}$, $j \geq 1$ are arbitrarily close 
to the weak 
unstable foliation. Thus we can replace the Jacobian of the flow by the weak 
unstable Jacobian, 
up to an $\eps>0$ error which can be taken as small as we like (possibly after increasing $T_0$ sufficiently). Thus
\[
J^{-T_0}_{\mc{L}_{jT_0}}(\Phi^{-T+jT_0}(\zbar)) 
\leq J^{wu}_{T_0}(\Phi^{-T+jT_0}(\zbar))(1+\eps)
\]
where $J^{wu}_t(m)$ is defined in \eqref{Ju}. But the right hand side is uniformly bounded by 
\[ \exp\Big(S_{T_0}(W_{\alpha_{j'}})\Big)(1 + \eps), \textrm{ with } \quad j':=T-jT_0.\]
Consequently, using \eqref{bJ;k}, \eqref{ccc} and \eqref{Jterms}, we find that \eqref{bbb} is bounded uniformly by
\[C (1+\eps)^{N}\exp\Big(\sum_{j=1}^{N}\demi S_{T_0}(W_{\alpha_{j'}})\Big)\]
for some subsequence $(\alpha_{j'})_{j'}\in A_1^{N}$ and some $C>0$ depending only on $T_0$. 
Now summing over all $\alpha$ in $\mc{A}_L$, we clearly obtain the bound \ah{This needs expanding}
\[\begin{gathered}
\sum_{\alpha\in \mc{A}_L}||\chi U(s_0)\Pi_{\alpha_L}U_{\alpha_{L-1}}\dots 
U_{\alpha_1}\chi ||_{L^1\to L^\infty}\leq Ch^{-d/2}(1+\eps)^N\sum_{\alpha\in \mc{A}_{L}}
\prod_{j'=1}^{N}e^{\demi S_{T_0}(W_{\alpha_j'})}\\
\leq Ch^{-d/2}(1+\eps)^N\Big(\sum_{a\in A_1}e^{\demi S_{T_0}(W_a)}\Big)^{N}
\end{gathered}\]
which from \eqref{epsilon} proves the Lemma since $NT_0$ is comparable to $T$. 
\qed\\

\textsl{Completion of the proof of Proposition~\ref{displogh}}:
We first note that for times $T \leq 1$, the estimate \ah{I added this first part because it seemed to me that in Lemma~\ref{mainest} we needed to have $NT_0$ comparable to $T$, i.e. $T \geq C T_0$.}
$$
||\chi \psi(h^2\Delta)e^{iTh\Delta}\chi||_{L^1\to L^\infty}\leq C(Th)^{-n/2}
$$
follows from the parametrix construction in \cite[Section 2.2]{BGT}. 

For times $1 \leq T \leq 4T_0$, the estimate can be obtained essentially as above: if $T=N+s_0$ with $N\in\nn$ and $s_0\in(0,1)$, the operator $\chi U(T)\psi(h^2\Delta)\chi$  is, modulo $O(h^\infty))$, a finite sum of terms of the form 
\begin{equation}\label{T<4T_0}
\chi U(T)\psi(h^2\Delta) \chi = \chi \psi(h^2\Delta)U(s_0)U_{\alpha_N}U_{\alpha_{N-1}}\dots U_{\alpha_1} \chi
\end{equation} 
where $U_{\alpha_j}=U(1)\Pi_{\alpha_j}\psi_+(h^2\Delta)$ and $\alpha_j\in A$ like above. Using the assumption that 
the region $\pi(W_0)$ is geodesically convex and that the support of $\chi$ is included in $M_-$, 
we see that only the cases where all the $\alpha_j$ are non-zero is not $O(h^\infty)$. But then, since the $\Pi_{\alpha_j}$ are microsupported in the part of the manifold which has negative curvature, then by the method of Anantharaman \cite{An} as we just explained before, the operators of 
\eqref{T<4T_0} are Lagrangian distributions and enjoy the $L^1\to L^\infty$ estimate
\[||\chi \psi(h^2\Delta)U(s_0)U_{\alpha_N}U_{\alpha_{N-1}}\dots U_{\alpha_1} \chi||_{L^1\to L^\infty}\leq C (Th)^{-\ndemi}.\]
and we sum those finitely many terms to obtain the desired result. 


For $T \geq 4 T_0$ we can apply  Lemmas \ref{lemrn}, \ref{lemAN} with the estimate of Lemma \ref{mainest} where 
$\eps_0+\eps$ is chosen smaller than $-P(1/2)$; we obtain the estimate 
\[||\chi \psi(h^2\Delta)e^{iTh\Delta}\chi||_{L^1\to L^\infty}\leq Ch^{-n/2}e^{-\beta T}\]
for some $\beta>0$, and all $T\in (0,\log(1/h))$. It suffices to set $t=Th$
and we get the desired result since $e^{-\beta T}T^{n/2}\leq C$.
\qed\\

\subsection{Proof of Theorem~\ref{sharpstrich}}\label{sec-proof} We shall be brief here since the proof was outlined already in Section~\ref{strategy}. We use the notation from that section. Thus, $u_j = \varphi(t/h|\log h| - j) \chi u$
satisfies 
$$
(i \partial_t - \Delta_M) u_j = w_j' + w_j''
$$
where $w_j',w_j''$ are defined in \eqref{decomegaj}.
Choose $\chi_- \in C_0^\infty(M)$ supported in $M_-$ and identically $1$ on the support of $\chi$, and $\chi_+ \in C^\infty(M)$ so that $1 - \chi_+ \in C_0^\infty(M)$ is identically $1$ on $\pi(K)$ and is $0$ on the support of $\nabla \chi$. Then $u_j = \chi_+u_j $ and $w_j' = \chi_- w_j'$, $w_j'' = \chi_+ w_j''$. We define $u_j'$ by 
 \begin{equation}
u_j'(t) = \chi_- \int_{(j-1)h|\log h|}^{t} e^{-i(t-s) \Delta_M} \chi_+ w_j'(s) \, ds 
\end{equation}
with $u_j''$ defined analogously. Clearly $u_j' + u_j'' = u_j$.

To treat $w_j''$, consider $\tilde u_j''$ defined by 
 \begin{equation}
\tilde u_j''(t) = \chi_- \int_{(j-1)h|\log h|}^{(j+1)h|\log h|} e^{-i(t-s) \Delta_M} \chi_+ w_j''(s) \, ds = \chi_-e^{-it\Delta_M} \int_{(j-1)h|\log h|}^{(j+1)h|\log h|} e^{is \Delta_M} \chi_+ w_j''(s) \, ds.
\end{equation}
Using Lemma~\ref{T} and Remark~\ref{remark-conic} we see that  
$$
\| \int_{(j-1)h|\log h|}^{(j+1)h|\log h|} e^{is \Delta_M} \chi_+ w_j(s) \, ds \|_{L^2(M)} \leq C h^{1/2} \| w_j'' \|_{L^2_t; L^2(M)}.
$$
Then Theorem~\ref{scsl} applied to this $L^2$ function shows that
\begin{equation}
\| \tilde u_j'' \|_{L^p_t; L^q(M)} \leq C h^{1/2} \| w_j'' \|_{L^2_t; L^2(M)}
\end{equation}
and the same estimate holds for $u_j''$ by Christ-Kiselev. 
To treat $w_j'$, consider $\tilde u_j'$ defined by 
 \begin{equation}
\tilde u_j'(t) = \chi_- \int_{(j-1)h|\log h|}^{(j+1)h|\log h|} e^{-i(t-s) \Delta_M} \chi_+ w_j'(s) \, ds = \chi_-e^{-it\Delta_M} \int_{(j-1)h|\log h|}^{(j+1)h|\log h|} e^{is \Delta_M} \chi_+ w_j'(s) \, ds.
\end{equation}
The dual estimate to Theorem~\ref{lsll} implies 
$$
\| \int_{(j-1)h|\log h|}^{(j+1)h|\log h|} e^{is \Delta_M} \chi_+ w_j'(s) \, ds \|_{L^2(M)} \leq 
C \big(h|\log h|\big)^{1/2} \| w_j' \|_{L^2_t; L^2(M)}\leq \frac{C\| \chi u\|_{L^2_t;L^2(M)}}{(h|\log h|)^{1/2}}.
$$
using also $\omega_j^{'}=i(h|\log h|)^{-1}\varphi'(t/h|\log h|-j)\chi u$.
Then  we can use Theorem~\ref{scsl} applied to this $L^2$ function shows that
\begin{equation}
\| \tilde u_j' \|_{L^p_t; L^q(M)} \leq C(h|\log h|)^{-1/2}  \| \chi u \|_{L^2_t; L^2(M)}
\end{equation}
and the same estimate holds for $u_j'$ by Christ-Kiselev. 

Squaring and summing over $j$ gives 
\begin{equation}
\sum_{j=1}^{N-1} \| u_j \|_{L^p_t; L^q(M)}^2 \leq 
C \sum_{j=1}^{N-1} \Big( h \| w_j'' \|_{L^2_t;L^2(M)}^2 + 
\frac{1}{h|\log h|} \| w_j' \|_{L^2_t;L^2(M)}^2 \Big)
\end{equation}
and the right hand side is no bigger than $C \| u_0 \|_{L^2(M)}^2$ using Lemma~\ref{T} for $w_j''$ and Theorem~\ref{lsll} for $w_j'$. 
Using the continuous embedding from $l^2(\NN)$ to $l^p(\NN)$ as in Section~\ref{strategy} gives
$$
\| \chi e^{-it\Delta_M} \psi(h^2 \Delta) u_0 \|_{L^p[0,1]; L^q(M)} \leq C \| u_0 \|_{L^2(M)}.
$$
Together with Theorem~\ref{htwt} this gives the Strichartz estimate without the space cutoff $\chi$:
$$
\| e^{-it\Delta_M} \psi(h^2 \Delta) u_0 \|_{L^p[0,1]; L^q(M)} \leq C \| u_0 \|_{L^2(M)}.
$$
Finally using Bouclet's Littlewood-Paley estimate (equation (1.4) of \cite{Bo}) and the argument in \cite{BGT}, we remove the frequency cutoff and obtain \eqref{sharpstrich3}, which completes the proof. 

\begin{rem}
The restriction $p > 2$ in Theorem~\ref{zzzz} is only required because we use the Christ-Kiselev lemma. It is likely that this condition could be eliminated (for $d > 2$) with a more careful analysis. 
\end{rem}

\end{document}